\documentclass[12pt]{article}
\usepackage{graphicx}
\usepackage{amsmath,amsthm,amssymb,enumerate}
\usepackage{euscript,mathrsfs}
\usepackage{color}
\usepackage{dsfont}
\usepackage{url}
\usepackage[left=2cm,right=2cm,top=3.5cm,bottom=3.5cm]{geometry}
\usepackage{color}
\usepackage[framemethod=tikz]{mdframed}
\allowdisplaybreaks
\usepackage{esint}

\usepackage{soul}

\catcode`\@=11 \@addtoreset{equation}{section}

\catcode`\@=12

\newtheorem{Theorem}{Theorem}[section]
\newtheorem{Proposition}[Theorem]{Proposition}
\newtheorem{Lemma}[Theorem]{Lemma}
\newtheorem{Corollary}[Theorem]{Corollary}

\theoremstyle{definition}
\newtheorem{Definition}[Theorem]{Definition}

\newtheorem{Remark}[Theorem]{Remark}

\newcommand{\bTheorem}[1]{
	\begin{Theorem} \label{T#1} }
	\newcommand{\eT}{\end{Theorem}}

\newcommand{\bProposition}[1]{
	\begin{Proposition} \label{P#1}}
	\newcommand{\eP}{\end{Proposition}}

\newcommand{\bLemma}[1]{
	\begin{Lemma} \label{L#1} }
	\newcommand{\eL}{\end{Lemma}}

\newcommand{\bCorollary}[1]{
	\begin{Corollary} \label{C#1} }
	\newcommand{\eC}{\end{Corollary}}

\newcommand{\bRemark}[1]{
	\begin{Remark} \label{R#1} }
	\newcommand{\eR}{\end{Remark}}

\newcommand{\bDefinition}[1]{
	\begin{Definition} \label{D#1} }
	\newcommand{\eD}{\end{Definition}}

\newcommand{\wvuh}{\widetilde{\vu}_h}
\newcommand{\wTheta}{\widetilde{\mathcal{T}}}

\newcommand{\Del}{\Delta_x}

\newcommand{\MTC}{\wtilde{\mathcal{T}}}
\newcommand{\vrB}{\vr_B}

\newcommand{\Ds}{\mathbb{D}_x}

\newcommand{\tvre}{\tvr_\ep}
\newcommand{\tvte}{\tvt_\ep}
\newcommand{\tvue}{\tvu_\ep}

\newcommand{\avintO}[1]{\fint_{\Omega} #1 \dx}

\newcommand{\data}{{\rm data}}

\newcommand{\bfphi}{\boldsymbol{\varphi}}

\newcommand{\bFormula}[1]{
	\begin{equation} \label{#1}}
	\newcommand{\eF}{\end{equation}}

\newcommand{\vuh}{\vu_h}

\newcommand{\Divh}{{\rm div}_h}
\newcommand{\Gradh}{\nabla_h}

\newcommand{\Ov}[1]{\overline{#1}}

\newcommand{\vr}{\varrho}
\newcommand{\vre}{\vr_\ep}

\newcommand{\vte}{\vt_\ep}
\newcommand{\vue}{\vu_\ep}
\newcommand{\tvr}{\wtilde \vr}
\newcommand{\tvu}{{\wtilde \vu}}
\newcommand{\tvt}{\wtilde \vt}

\newcommand{\vt}{\vartheta}
\newcommand{\vu}{\vc{u}}

\newcommand{\vc}[1]{{\bf #1}}

\newcommand{\Div}{{\rm div}_x}
\newcommand{\Grad}{\nabla_x}

\newcommand{\dx}{\,{\rm d} {x}}

\newcommand{\dt}{\,{\rm d} t }

\newcommand{\vU}{\vc{U}}

\newcommand{\intO}[1]{\int_{\Omega} #1 \ \dx}

\newcommand{\intOh}[1]{\int_{B(r)} #1 \ \dx}

\newcommand{\D}{{\rm d}}

\newcommand{\ep}{\varepsilon}

\newcommand{\vtB}{\vt_B}

\newcommand{\br}{ \nonumber \\ }

\def\softd{{\leavevmode\setbox1=\hbox{d}%
		\hbox to 1.05\wd1{d\kern-0.4ex{\char039}\hss}}}
\definecolor{Cgrey}{rgb}{0.85,0.85,0.85}
\definecolor{Cblue}{rgb}{0.50,0.85,0.85}
\definecolor{Cred}{rgb}{1,0,0}
\definecolor{fancy}{rgb}{0.10,0.85,0.10}
\definecolor{amaranth}{rgb}{0.9, 0.17, 0.31}

\newcommand\Cbox[2]{%
	\newbox\contentbox%
	\newbox\bkgdbox%
	\setbox\contentbox\hbox to \hsize{%
		\vtop{
			\kern\columnsep
			\hbox to \hsize{%
				\kern\columnsep%
				\advance\hsize by -2\columnsep%
				\setlength{\textwidth}{\hsize}%
				\vbox{
					\parskip=\baselineskip
					\parindent=0bp
					#2
				}%
				\kern\columnsep%
			}%
			\kern\columnsep%
		}%
	}%
	\setbox\bkgdbox\vbox{
		\color{#1}
		\hrule width  \wd\contentbox %
		height \ht\contentbox %
		depth  \dp\contentbox
		\color{black}
	}%
	\wd\bkgdbox=0bp%
	\vbox{\hbox to \hsize{\box\bkgdbox\box\contentbox}}%
	\vskip\baselineskip%
}

\mdfdefinestyle{MyFrame}{%
	linecolor=black,
	outerlinewidth=1pt,
	roundcorner=5pt,
	innertopmargin=\baselineskip,
	innerbottommargin=\baselineskip,
	innerrightmargin=10pt,
	innerleftmargin=10pt,
	backgroundcolor=white!20!white}

\newcommand{\veps}{\varepsilon}

\newcommand{\wtilde}{\widetilde}

\newcommand{\lan}{\langle}
\newcommand{\ran}{\rangle}

\newcommand{\aleq}{\lesssim}

\allowdisplaybreaks

\begin{document}


\title{\bf Continuous data assimilation applied to 
\\
the Rayleigh--B\' enard problem 
\\
for compressible fluid flows}

\author{Eduard Feireisl	
\thanks{The work of E.F. was supported by the
		Czech Sciences Foundation (GA\v CR), Grant Agreement 24-11034S.
		The Institute of Mathematics of the Academy of Sciences of
		the Czech Republic is supported by RVO:67985840.} 
and Wladimir Neves
\thanks{W. Neves has received research grants from CNPq
	through the grants 313005/2023-0, 406460/2023-0, and also by FAPERJ 
	(Cientista do Nosso Estado) through the grant E-26/204.171/2024
	}
}

\date{}

\maketitle

\medskip

\centerline{Institute of Mathematics of the Academy of Sciences of the Czech Republic}

\centerline{\v Zitn\' a 25, CZ-115 67 Praha 1, Czech Republic}

\bigskip

\centerline{Instituto de Matem\' atica, Universidade Federal do Rio de Janeiro}

\centerline{C.P. 68530, Cidade Universit\' aria 21945-970, Rio de Janeiro, Brazil}

\maketitle

\begin{abstract}
We apply a continuous data assimilation method to the Navier--Stokes--Fourier 
system governing the evolution of a compressible, rotating and thermally driven fluid. 
A rigorous proof of the tracking property is given in the asymptotic regime of low Mach 
and high Rossby and Froude numbers. Large data in the framework of weak solutions are considered. 
\end{abstract}


{\small

\noindent

\medbreak


}

\section{Introduction}
\label{i}

Continuous Data Assimilation is a framework to incorporate observational data continuously into 
a mathematical model represented usually by an evolutionary system of partial differential equations. 
The method is often based on adding \emph{nudging} terms determined by interpolated observations to the governing system of equations. 
Although the leading idea is conceptually simple and has been described in the seminal paper by Hoke and Anthes 
\cite{HoAn} as early as in 1976, the first rigorous mathematical treatment appeared much later in the pioneering work 
of Azouani, Olson, and Titi \cite{AzOlTi}. 

\medskip
Motivated by the existing theory of determining modes and 
finite dimensional attractors for the Navier--Stokes system (Bercovici et al. \cite{BeCoFoMa}, Constantin, Foias and Temam \cite{CFT}, 
Constantin et al. \cite{CoFoMaTe}), they proposed a nudging term based on interpolants ranging 
in a suitable finite dimensional space. Applicability of the method to the real world problems arising e.g. in meteorology is, however, 
limited by the lack of the required well--posedness theory for the physically relevant 3-D Navier--Stokes system. 
The latter problem has been recently attacked in a series of works of Biswas et al. \cite{BalBis2}, \cite{BalBis}, \cite{BisPri}, where 
a new regularity criterion for the 3-D Navier--Stokes system is proposed based solely on the observed values of the 
interpolation functionals. Note, however, that such a criterion guarantees regularity of the 
background solution only in the time interval when the observation 
data are collected, but not after the prediction is made or computed.   

\subsection{Problem formulation}

The real world applications arising in meteorology, astrophysics or geophysics are often represented by 
the general Navier-Stokes-Fourier system describing the motion of a compressible, viscous and heat conducting fluid confined 
to a physical domain and supplemented by the relevant boundary conditions. We consider a rotating fluid system  
confined to a bounded domain $\Omega \subset R^3$, 
\begin{equation} \label{s1}
\Omega = B(r) \times (0,1),\qquad  B(r) = \left\{ {x}_h = (x_1, x_2) \ \Big| \ |x_h | < r \right\}.  
\end{equation}
Here and hereafter, we write $x \in\Omega$ in the form $x = ({x}_h, x_3)$ to stress the anisotropy between ``horizontal'' and ``vertical'' variables, ${x}_h$ and $x_3$, respectively.  
We assume the fluid domain is rotating around its vertical axis. 

Let $\vr = \vr(t,x)$, $\vt = \vt(t,x)$, and $\vu = \vu(t,x)$ denote the mass density, the (absolute) temperature, and the velocity field, 
respectively.	
The classical principles of conservation of mass, linear momentum, and energy written in the rotating frame give rise to the 
\emph{Navier--Stokes--Fourier} (NSF) system: 
\begin{align} 
	\partial_t \vr + \Div (\vr \vu) &= 0, \label{s2}\\
	\partial_t (\vr \vu) + \Div (\vr \vu \otimes \vu) + \frac{1}{\ep^2} \Grad p(\vr, \vt) + \frac{1}{\sqrt{\veps}} \vc{e}_3\times \vr \vu
	&= \Div \mathbb{S}(\vt, \Ds \vu) + \frac{1}{\ep} \vr \Big( \Grad G  +  \Grad |{x}_h |^2 \Big) , \label{s3} \\
	\partial_t (\vr e(\vr, \vt)) + \Div (\vr e(\vr, \vt) \vu) + \Div  \vc{q} (\vt, \Grad \vt) &= 
	\ep^2 \mathbb{S} (\vt, \Ds \vu) : \Ds \vu - p(\vr, \vt) \Div \vu .
	\label{s4}	
\end{align}

\noindent
The symbols $p = p(\vr, \vt)$ and $e = e(\vr, \vt)$ denote the pressure and the internal energy, respectively, specified constitutively in Section \ref{m} below. 
The function $G$
represents the gravitational potential acting in the vertical direction,
\begin{equation} \label{s5}
	G = - x_3 .
\end{equation}
The effect of rotation is represented by the Coriolis force $\vc{e}_3\times \vr\vu$, where $\vc{e}_3=(0,0,1)$, and the centrifugal 
force  $\vr \Grad|{x}_h|^2$.
The viscous stress tensor is given by \emph{Newton's rheological law} 
\begin{equation} \label{s6}
	\mathbb{S}(\vt, \Ds \vu) = 2\mu(\vt) \left( \Ds\vu - \frac{1}{3} \Div \vu \mathbb{I} \right) + \eta(\vt) \Div \vu \mathbb{I}, 
\end{equation}
where $\Ds\vu = \frac{1}{2}(\Grad \vu + \Grad^t \vu)$ is the symmetric part of the velocity gradient $\Grad \vu$. 
Analogously, the heat flux is given by \emph{Fourier's law}
\begin{equation} \label{s7} 
	\vc{q}(\vt, \Grad \vt) = - \kappa (\vt) \Grad \vt.
\end{equation}
The system \eqref{s2}--\eqref{s4} is scaled by a small parameter $\ep > 0$, specifically, the fluid flow is nearly 
incompressible (the Mach number proportional to $\ep$), stratified (the Froude number $\sqrt{\ep}$), and rotating 
(the Rossby number ${\ep}^{-\frac{1}{2}}$ enforcing the centrifugal force proportional $\ep^{-1}$), see \cite{FanFei2025}
for the physical background.

We further suppose that the fluid motion is driven by thermal convection. In accordance with the scaling, we impose 
the Rayleigh--B\' enard type boundary condition 
\begin{equation} \label{s9} 
\vt|_{\partial \Omega} = \Ov{\vt} + \ep \vtB, 
\end{equation}
where $\Ov{\vt} > 0$ is a positive constant representing the background temperature. In contrast with $\Ov{\vt}$, the temperature 
profile $\vtB = \vtB(x)$ need not be positive as the scaling parameter $\ep$ vanishes in the asymptotic regime. In addition, we impose 
the no--slip boundary conditions for the velocity on the lateral boundary 
\begin{equation} \label{s10}
\vu|_{ \partial B(r) \times (0,1) } = 0, 
\end{equation}
and, to avoid the friction effect of a boundary layer, the complete slip boundary conditions 
\begin{equation} \label{s11}
\vu \cdot \vc{n}|_{ B(r) \times \{ x_3 = 0, 1 \} } = 0,\ 
[\mathbb{S}(\vt, \Ds \vu) \cdot \vc{n} ] \times \vc{n} |_{ B(r) \times \{ x_3 = 0, 1 \} } = 0
\end{equation} 
on the horizontal parts of the boundary. 

The NSF system \eqref{s2}--\eqref{s11} will be termed \emph{observed system} and its solution \emph{observed solution}.
The time evolution of the observed system is considered on the time interval 
\begin{equation} \label{s8}
t \in (T^-, T^+),\ \mbox{where}\  T^- < 0,\ T^+ > 0. 
\end{equation}
The initial state (data) of the system is {\it a priori} not prescribed, however, we suppose the  
solution emanates from the state that is \emph{well prepared} with respect to the adopted scaling. Specifically, 
there is a reference constant density $\Ov{\vr} > 0$ such that   
\begin{align}
\left\| \frac{\vr( T^-, \cdot) - \Ov{\vr} }{\ep} \right\|_{L^\infty(\Omega)} \aleq 1 ,\  
\left\| \frac{\vr( T^-, \cdot) - \Ov{\vr}}{\ep} -  {\mathcal{R}}_0 \right\|_{L^1(\Omega)} \leq d, \br 
\left\| \frac{\vt( T^-, \cdot) - \Ov{\vt}}{\ep} \right\|_{L^\infty(\Omega)} \aleq 1 ,\  
\left\| \frac{\vt( T^-, \cdot) - \Ov{\vt}}{\ep} -  \mathcal{T}_0 \right\|_{L^1(\Omega)} \leq d, \br
\left\| \vu( T^-, \cdot) \right\|_{L^\infty(\Omega; R^3)} \aleq 1,\ 
\left\| \vu( T^-, \cdot) - {\vu}_{0,h} \right\|_{L^1(\Omega; R^3)} \leq d
\label{s12}
\end{align}
for certain functions $
{\mathcal{R}}_0$, $\mathcal{T}_0$, and $\vu_{0,h}$ satisfying
\begin{align} 
\mathcal{T}_0 &\in W^{2,\infty}(\Omega), \ {\mathcal{T}}_0|_{\partial \Omega} = \vtB, \br 
[{\vu}_{0,h}, 0]\equiv {\vu}_{0,h} 
&\in W^{2,\infty}(B(r); R^2),\  
{\vu}_{0,h}|_{\partial B(r)} = 0,\ \Divh {\vu}_{0,h} = 0, 	
\label{s13}
\end{align}
and	
\begin{equation} \label{s14}
\frac{\partial p(\Ov{\vr}, \Ov{\vt})}{\partial \vr } \Grad {\mathcal{R}}_0 + 
\frac{\partial p(\Ov{\vr}, \Ov{\vt})}{\partial \vt } \Grad {\mathcal{T}}_0 = \Ov{\vr} \Grad \left( G +  |x_h|^2 \right).	
\end{equation}

\begin{Remark} \label{iR1}
Here and throughout the paper, the symbol $a\leq b$ means there exists a positive constant $C$ such that $a \leq Cb$. The subscript $h$ 
appended to a differential operator means the latter acts only in the horizontal variable $x_h$. 
\end{Remark}	
	
\noindent	
Identity \eqref{s14} is nothing other than the celebrated Boussinesq relation arising from the linearization of 
the momentum equation around the reference state $(\Ov{\vr}, \Ov{\vt})$ 
in the limit regime $\ep \to 0$, see e.g. Zeytounian \cite{ZEY1}.
It is worth--noting that the well--prepared velocity field $\vu_{0,h}$ is planar depending solely on the horizontal variable 
$x_h$ as a consequence of the fast rotation, see \cite{FanFei2025} for a detailed discussion of this property that is somehow 
at odds with the commonly accepted scenario advocated in applications (see e.g. Ecke and Shiskhina \cite{EckShi}).

\subsection{Continuous data assimilation problem}

The leading idea of Continuous Data Assimilation is that the data are generated by an exact solution of $(\vr, \vt, \vu)$ of  
the observed NSF system considered on a time interval $(T^-,T^+)$. The data are collected and implemented  
in a \emph{synchronized system} in a time interval $[0, T)$,
\[
T^- < 0 < T < T^+.  
\]
The synchronized system then generates a solution $(\tvr, \tvt, \tvu)$ in $[0, T^+)$ \emph{predicting} the 
expected behavior of the observed solution in the time interval $[T, T^+)$.
 
The \emph{synchronized} solution $(\tvr, \tvt, \tvu)$ satisfies the NSF system of field equations 
augmented by nudging terms of ``friction'' type active on the synchronization time interval $(0,T)$, see Section \ref{m}.
Motivated by the seminal paper by Azouani, Olson and Titi \cite{AzOlTi}, we consider the nudging terms in the form
\[
- \Lambda \left( I_\delta [\tvu] - I_\delta [\vu] \right) \mathds{1}_{t \in (0,T)}
\]
acting as a friction force in the momentum equation \eqref{s3}, and
\[
- \tvt  \Lambda \left( I_\delta [\tvt] - I_\delta [\vt] \right) \mathds{1}_{t \in (0,T)}
\] 
as a heat source in the internal energy balance \eqref{s4}. 
Here $\Lambda > 0$ denotes the \emph{nudging parameter}, whereas $I_\delta$ are the interpolants representing the available ``measurements'' 
of the  observed solution $\vu$, $\vt$. We focus on linear interpolants 
that can be represented as \emph{projections} on a suitable finite--dimensional space approaching identity 
for $\delta \to 0$.  
Generalization 
of our results to non--linear interpolants considered e.g. by Carlson, Larios, and Titi \cite{CaLaTi} is possible  
under adequate structural restrictions.

Our main result, formulated in Section \ref{m}, Theorem \ref{mT1}, asserts the existence of $\Lambda > 0$, 
$\delta > 0$, and  $\ep_0 > 0$, $d_0 > 0$ such that the synchronized solution $(\tvr, \tvt, \tvu)$ is arbitrarily close to 
the observed solution $(\vr, \vt, \vu)$ on the time interval $(T, T^+)$ as long as 
$0 < \ep < \ep_0$, $0 < d < d_0$. To the best of our knowledge, this is the first proof of convergence
of a continuous data assimilation method applied to the physically realistic 3-D compressible NSF system in the Rayleigh--B\' enard 
convection regime. In particular, the results is stated in the framework of weak solutions, where the desired well--posedness 
remains an outstanding open problem. 

Our strategy leans on several steps: 

\begin{itemize}
\item Using \cite{FanFei2025} we identify the asymptotic limit - \emph{the target problem} - of the observed NSF system for $\ep \to 0$ in the time interval $(T^-,T^+)$.

\item Adopting the standard energy estimates to the target problem, we show the convergence of the continuous data assimilation method for the latter.

\item Modifying the arguments of \cite{FanFei2025} we identify the asymptotic limit of the synchronized NSF system in the time interval $(0,T^+)$. In particular, we show solutions of the synchronized NSF system and the synchronized system 
associated to the target problem are close in $[0, T^+)$.

\item Summing up the previous results, we show proximity of the observed and synchronized solutions in the time interval $(T, T^+)$.
\end{itemize}

The paper is organized as follows. In Section \ref{m}, we collect the preliminary results concerning the theory of weak 
solutions to the NSF system and state the main result. In Section \ref{as}, using the results of \cite{FanFei2025}, we identify the asymptotic limit of solutions 
to the observed system. In Section \ref{das}, we apply the data assimilation method proposed by Azouani, Olson, Titi \cite{AzOlTi} to the target problem. Finally, in Section \ref{F}, we perform the asymptotic limit for the synchronized NSF system.
 
\section{Preliminaries, main results}
\label{m}

We collect the necessary piece of information concerning the theory of weak solutions of the NSF system and state our main result.

\subsection{Hypotheses imposed on constitutive relations}
\label{HIC}

A suitable form of constitutive relations was introduced in \cite[Chapters 1,2]{FeNo6A} and \cite[Chapter 1]{FeiNovOpen}.
We refer the reader to \cite[Chapter 1]{FeNo6A} for their physical background. 

The pressure equation of state (EOS) reads
\[
p(\vr, \vt) = p_{\rm m} (\vr, \vt) + p_{\rm rad}(\vt), 
\]
where $p_{\rm m}$ is the pressure of a general \emph{monoatomic} gas, 
\begin{equation} \label{con1}
	p_{\rm m} (\vr, \vt) = \frac{2}{3} \vr e_{\rm m}(\vr, \vt),
\end{equation}
supplemented with the radiation pressure 
\[
p_{\rm rad}(\vt) = \frac{a}{3} \vt^4,\qquad a > 0.
\]
The pressure and the internal energy are interrelated through Gibbs' law
\begin{equation} \label{Gibbs}
	\vt Ds = De + p D \left( \frac{1}{\vr} \right),
\end{equation}
{where $s$ is the (specific) entropy.} {In addition, we impose the hypothesis of thermodynamic stability}
\begin{equation} \label{ThSt}
	\frac{\partial p(\vr, \vt)}{\partial \vr} > 0,\ \frac{\partial e(\vr, \vt)}{\partial \vt} > 0. 
\end{equation} 

The main structural hypotheses concerning EOS are formulated below:

\begin{itemize}
	
	\item {Gibbs' relation} together with \eqref{con1} yield the specific form of the pressure $p_m$, 
	\[
	p_{\rm m} (\vr, \vt) = \vt^{\frac{5}{2}} P \left( \frac{\vr}{\vt^{\frac{3}{2}}  } \right)
	\]
	for a certain $P \in C^1[0,\infty)$.
	Consequently, we get
	\begin{equation} \label{w9}
		p(\vr, \vt) = \vt^{\frac{5}{2}} P \left( \frac{\vr}{\vt^{\frac{3}{2}}  } \right) + \frac{a}{3} \vt^4,\quad
		e(\vr, \vt) = \frac{3}{2} \frac{\vt^{\frac{5}{2}} }{\vr} P \left( \frac{\vr}{\vt^{\frac{3}{2}}  } \right) + \frac{a}{\vr} \vt^4, \qquad a > 0.
	\end{equation}
	
	\item Hypothesis of thermodynamic stability \eqref{ThSt}
	expressed in terms of  $P$ yields
	\begin{equation} \label{w10}
		P \in C^1[0, \infty),\ P(0) = 0, \ P'(Z) > 0 \ \mbox{for}\ Z \geq 0,\qquad 0<  \frac{ \frac{5}{3} P(Z) - P'(Z) Z }{Z} \leq c \ \mbox{ for }\ Z > 0,
	\end{equation} 	
	where the interpretation of the upper bound is the boundedness of the specific heat at constant volume.	
	\item 
	The specific {entropy} takes the form 
	\begin{equation} \label{w12}
		s(\vr, \vt) = s_{\rm m}(\vr, \vt) + s_{\rm rad}(\vr, \vt),\qquad s_{\rm m} (\vr, \vt) = \mathcal{S} \left( \frac{\vr}{\vt^{\frac{3}{2}} } \right),\qquad
		s_{\rm rad}(\vr, \vt) = \frac{4a}{3} \frac{\vt^3}{\vr}, 
	\end{equation}
	where 
	\begin{equation} \label{w13}
		\mathcal{S}'(Z) = -\frac{3}{2} \frac{ \frac{5}{3} P(Z) - P'(Z) Z }{Z^2} < 0.
	\end{equation}
	
	\item 
	Finally, 
	we impose  requiring  the entropy to vanish 
	when the absolute temperature approaches zero, 
	\begin{equation} \label{w14}
		\lim_{Z \to \infty} \mathcal{S}(Z) = 0.
	\end{equation}
	In addition, we suppose 
	\begin{equation} \label{w14a}
		P \in C^1[0,\infty) \ \mbox{is such that } \ \liminf_{Z\to\infty}\frac{P(Z)}{Z}>0, 
	\end{equation}
	see \cite[Section 2.1.1]{FeiLuSun} for details.
	
\end{itemize}

\begin{Remark} \label{crR}
	
It is worth noting that the above set of hypotheses imposed on 
$p_m$, with the exception of 
\eqref{w14}, is compatible with the standard Boyle-Mariotte EOS  
$p_m(\vr, \vt) = \vr \vt$.	
	
\end{Remark}	

As for the transport coefficients, we suppose they are continuously differentiable functions of the 
temperature satisfying
\begin{align} 
	0 < \underline{\mu}(1 + \vt) &\leq \mu(\vt),\qquad |\mu'(\vt)| \leq \Ov{\mu}, 
	\br 
	0 &\leq \underline{\eta} (1 + \vt) \leq \eta (\vt) \leq \Ov{\eta}(1 + \vt), 
	\br
	0 < \underline{\kappa} (1 + \vt^\beta) &\leq \kappa (\vt) \leq \Ov{\kappa}(1 + \vt^\beta), 
	\quad \mbox{ where }\ \beta > 6. \label{w16}
\end{align}

\subsection{Weak solutions to the observed system}

Having collected all necessary hypotheses, we are ready to introduce the concept of \emph{weak solution} to the observed NSF system.

\begin{Definition}[{\bf Weak solution to the NSF system}] \label{DL1}
	We say that a trio $(\vr, \vt, \vu)$ is a \emph{weak solution} of the observed NSF system \eqref{s2}--\eqref{s7},
	accompanied with the boundary conditions \eqref{s9}--\eqref{s11}, and
	the initial data
	\[
	\vr(T^-, \cdot) = \vr_{0},\ \vt(T^-, \cdot) = \vt_{0},\ 
		\vu(T^-, \cdot) = \vu_{0}
	\]
	if the following holds:
	\begin{itemize}
		
		\item The solution belongs to the {regularity class} 
		\begin{align}
			\vr &\in L^\infty(T^-,T^+; L^{\frac{5}{3}}(\Omega)),\ \vr \geq 0 
			\ \mbox{a.a.~in}\ (T^-,T^+) \times \Omega, \br
			\vu &\in L^2(T^-,T^+; W^{1,2} (\Omega; R^3)), \ \vu|_{\partial B(r)} = 0,\ \vu \cdot \vc{n}|_{x_3 = 0,1} = 0, \br 
			\vt^{\beta/2} ,\ \log(\vt) &\in L^2(T^-,T^+; W^{1,2}(\Omega)) \ \mbox{for some}\ \beta \geq 2,\ 
			\vt > 0 \ \mbox{a.a.~in}\ (T^-,T^+) \times \Omega, \br
			\vt&\in L^2(T^-,T^+; W^{1,2} (\Omega)),\ \vt|_{\partial \Omega} = \Ov{\vt} + \ep \vtB.
			\label{m1}
		\end{align}
		
		\item The {equation of continuity} \eqref{s2} is satisfied in the sense of distributions,
		\begin{align} 
			\int_{T^-}^{T^+} \intO{ \Big[ \vr \partial_t \varphi + \vr \vu \cdot \Grad \varphi \Big] } \dt &=  - 
			\intO{ \vr_{0} \varphi(T^-, \cdot) }
			\label{m2}
		\end{align}
		for any $\varphi \in C^1_c([T^-,T^+) \times \Ov{\Omega} )$.
		\item The {momentum equation} \eqref{s3} is satisfied in the sense of distributions, 
		\begin{align}
			\int_{T^-}^{T^+} &\intO{ \left[ \vr \vu \cdot \partial_t \bfphi + \vr \vu \otimes \vu : \Grad \bfphi - {\frac{1}{\sqrt{\ep}} (\bf{e}_3 \times \vr \vu) \cdot 
					\bfphi} + 
				\frac{1}{\ep^2} p(\vr, \vt) \Div \bfphi \right] } \dt \br &= \int_{T^-}^{T^+} \intO{ \left[ \mathbb{S}(\vt, \Ds \vu) : \Grad \bfphi - \frac{1}{\ep} \vr \Grad G \cdot \bfphi - { \frac{1}{\ep} \vr \Grad |\vc{x}_h |^2 \cdot \bfphi }  \right] } \dt \br &- 
			\intO{ \vr_{0} \vu_{0} \cdot \bfphi (T^-, \cdot) }
			\label{m3}
		\end{align}	
		for any $\bfphi \in C^1_c([T^-, T^+) \times \Ov{\Omega}; R^3)$ such that $\bfphi|_{\partial B(r)} = 0$, 
		$\bfphi \cdot \vc{n}|_{x_3 = 0,1} = 0$.
		
		\item There is a non--negative Radon measure - the entropy production rate -  $\sigma \in \mathcal{M}^+([T^-,T^+] \times \Ov{\Omega})$, 
		\begin{equation} \label{m4}
			\sigma \geq \frac{1}{\vt} \left( \ep^2 \mathbb{S}(\vt, \Ds \vu) : \Ds \vu - \frac{\vc{q}(\vt ,\Grad \vt) \cdot \Grad \vt}{\vt}\right).
		\end{equation}
		The {entropy balance} is satisfied in the following sense: 
		\begin{align}
			- \int_{T^-}^{T^+} &\intO{ \left[ \vr s(\vr, \vt) \partial_t \varphi + \vr s (\vr,\vt) \vu \cdot \Grad \varphi + \frac{\vc{q} (\vt, \Grad \vt )}{\vt} \cdot 
				\Grad \varphi \right] } \dt \br &= \int_{T^-}^{T^+} \int_{\Omega}{ \varphi \ \D \sigma(t,x)}  + \intO{ \vr_{0} s(\vr_{0}, \vt_{0}) 
				\varphi (T^-, \cdot) } 
			\label{m5} 
		\end{align}
		for any $\varphi \in C^1_c([T^-, T^+) \times \Omega)$.

		\item  The {ballistic energy inequality}
		\begin{align}  
			- &\int_{T^-}^{T^+} \partial_t \psi	\intO{ \left[ \ep^2 \frac{1}{2} \vr |\vu|^2 + \vr e(\vr, \vt) - \Theta \vr s(\vr, \vt) \right] } \dt  \br 
			&+ \int_{T^-}^{T^+} \int_{\Ov{\Omega}} \psi \Theta\  \D \sigma (t,x) 
			\br
			&\leq 
			\int_{T^-}^{T^+} \psi \intO{ \left[ \ep \vr \vu \cdot \Grad G + { \ep \vr \vu \cdot \Grad |{x}_h |^2 } \right] } \dt \br  &- 
			\int_{T^-}^{T^+} \psi \intO{ \left[ 
				\vr s(\vr, \vt) \partial_t \Theta + \vr s(\vr, \vt) \vu \cdot \Grad \Theta  + \frac{\vc{q}(\vt, \Grad \vt)}{\vt} \cdot \Grad \Theta \right] } \dt \br 
			&+ \psi(T^-) \intO{ \left[ \frac{1}{2} \ep^2 \vr_{0} |\vu_{0}|^2 + \vr_{0} e(\vr_{0}, \vt_{0}) - \Theta (T^-, \cdot) \vr_{0} s(\vr_{0}, \vt_{0}) \right] }
			\label{m6}
		\end{align}
		holds true for any $\psi \in C^1_c [T^-, T^+)$, $\psi \geq 0$, and \emph{any} smooth extension $\Theta$ of the boundary temperature, 
		\[
		\Theta > 0 \ \mbox{in}\ [T^-,T^+] \times \Ov{\Omega},\ \Theta|_{\partial \Omega} = \Ov{\vt} + \ep \vtB.                          
		\]
	\end{itemize}
	
\end{Definition}

The concept of weak solution specified in Definition \ref{DL1} has been introduced in \cite{ChauFei}, \cite[Chapter 12]{FeiNovOpen}, where 
global--in--time existence is established 
under the hypotheses of Section \ref{HIC}.

\subsection{Weak solutions to the (forced) synchronized system}

In contrast with the observed NSF system, the synchronized NSF system is considered on the time interval $(0,T^+)$. The momentum as well as the internal 
energy balance are augmented by nudging forces active in the time lap $(0,T)$, $0 < T < T^+$. Accordingly, the classical formulation of the synchronized NSF system reads:
 \begin{align} 
	\partial_t \tvr + \Div (\tvr \tvu) &= 0, \label{m7}\\
	\partial_t (\tvr \tvu) + \Div (\tvr \tvu \otimes \tvu) + \frac{1}{\ep^2} \Grad p(\tvr, \tvt) + \frac{1}{\sqrt{\veps}} \vc{e}_3\times \tvr \tvu
	&= \Div \mathbb{S}(\tvt, \Ds \tvu) + \frac{1}{\ep} \tvr \Big( \Grad G  +  \Grad |{x}_h |^2 \Big) \br 
	&- \Lambda \Big( I_\delta [\tvu] - I_{\delta}[\vu] \Big) \mathds{1}_{t \in [0,T]}, \label{m8} \\ 
	\partial_t (\tvr e(\tvr, \tvt)) + \Div (\tvr e(\tvr, \tvt) \tvu) + \Div  \vc{q} (\tvt, \Grad \tvt) &= 
	\ep^2 \mathbb{S} (\tvt, \Ds \tvu) : \Ds \tvu - p(\tvr, \tvt) \Div \tvu \br
	&- \tvt \Lambda \Big( I_\delta [\tvt] - I_{\delta}[\vt] \Big) \mathds{1}_{t \in [0,T]}.
	\label{m9}	
\end{align}

The symbol $I_\delta$ denotes the interpolation operators forcing the synchronized solution $(\tvr, \tvt, \tvu)$ to be close to the observed solution  
$(\vr, \vt, \vu)$. Motivated by Azouani, Olson and Titi \cite{AzOlTi}, we suppose
\begin{align}
I_\delta: L^2(\Omega) \to L^2 \cap L^\infty (\Omega) \ \mbox{is an $L^2-$orthogonal projection}, \br 
\| I_\delta (v) - v \|_{L^2(\Omega)} \aleq \delta \| v \|_{W^{1,2}(\Omega)}. 
\label{m10}	
\end{align}	
We remark that, a non--linear variant in the spirit of Carlson, Larios, and Titi \cite{CaLaTi} is possible pending some restrictions that will become clear below. 
Moreover, we do not distinguished between the scalar form of $I_\delta$ 
in \eqref{m9}
and its vector valued variant applied component wise to the velocity field 
in \eqref{m8}.

The adequate boundary conditions are 
\begin{align} 
	\tvt|_{\partial \Omega} &= \Ov{\vt} + \ep \vtB, \label{c5} \\
	\tvu|_{ \partial B(r) \times (0,1) } &= 0, \label{c6} \\ 
	\tvu \cdot \vc{n}|_{ B(r) \times \{ x_3 = 0, 1 \} } &= 0,\ 
	[\mathbb{S}(\tvt, \Ds \tvu) \cdot \vc{n} ] \times \vc{n} |_{ B(r) \times \{ x_3 = 0, 1 \} } = 0. 
	\label{m11}
\end{align}

The initial conditions specified at the time $t= 0$ are as follows
\begin{align}
	\left\| \frac{\tvr( 0, \cdot) - \Ov{\vr} }{\ep} \right\|_{L^\infty(\Omega)} \aleq 1 ,\  
	\left\| \frac{\tvr( 0, \cdot) - \Ov{\vr}}{\ep} - \wtilde{\mathcal{R}}_0 \right\|_{L^1(\Omega)} \aleq  d, \br 
	\left\| \frac{\tvt( 0, \cdot) - \Ov{\vt} }{\ep}  \right\|_{L^\infty(\Omega)} \aleq 1 ,\  
	\left\| \frac{ \tvt( 0, \cdot) - \Ov{\vt}}{\ep} -  \wTheta_0 \right\|_{L^1(\Omega)} \aleq d, \br
	\left\| \tvu( 0, \cdot) \right\|_{L^\infty(\Omega; R^3)} \aleq 1,\ 
	\left\| \tvu( 0, \cdot) - {\tvu}_{0,h} \right\|_{L^1(\Omega; R^3)} \aleq d.
	\label{m12}
\end{align}
Moreover, 
\begin{align} 
	\wtilde{\mathcal{T}}_0 &\in W^{2,\infty}(\Omega), \ \wtilde{\mathcal{T}}_0|_{\partial \Omega} = \vtB, \br 
	[\tvu_{0,h}, 0] \equiv{\tvu}_{0,h} 
	  &\in W^{2,\infty}(B(r); R^2),\   
	\tvu_{0,h}|_{\partial B(r)} = 0,\ \Divh \tvu_{0,h} = 0, 	
	\label{m13}
\end{align}
and	
\begin{equation} \label{m14}
	\frac{\partial p(\Ov{\vr}, \Ov{\vt})}{\partial \vr } \Grad \wtilde{\mathcal{R}}_0 + 
	\frac{\partial p(\Ov{\vr}, \Ov{\vt})}{\partial \vt } \Grad \wTheta_0 = \Ov{\vr} \Grad \left( G +  |x_h|^2 \right).	
\end{equation}	

\begin{Remark} \label{cR1}
	
	As a matter of fact, the initial data for the synchronized problem can be conveniently \emph{chosen} to simplify eventual computations. For instance, we may consider
	\begin{align}
		\tvr(0, \cdot) &= \Ov{\vr} + \ep \wtilde{\mathcal{R}}_0, \ \tvt( 0, \cdot) = \Ov{\vt} + \ep \wTheta_0, \ \tvu(0, \cdot) = 0, \br
		\wTheta_0 &\in W^{2,\infty}(\Omega), \ \wtilde{\mathcal{T}}_0|_{\partial \Omega} = \vtB, \br
		\frac{\partial p(\Ov{\vr}, \Ov{\vt})}{\partial \vr } \Grad \wtilde{\mathcal{R}}_0 + \frac{\partial p(\Ov{\vr}, \Ov{\vt})}{\partial \vt } \Grad \wtilde{\mathcal{T}}_0 &= \Ov{\vr} \Grad \left( G +  |x_h|^2 \right).
		\label{m15}
	\end{align}
	The density perturbation $\wtilde{\mathcal{R}}_0$ can be determined in terms of a suitable extension of the boundary temperature 
	$\vtB$. Indeed we have 
	\begin{equation} \label{m16}
		\frac{\partial p(\Ov{\vr}, \Ov{\vt})}{\partial \vr } \wtilde{\mathcal{R}}_0 + \frac{\partial p(\Ov{\vr}, \Ov{\vt})}{\partial \vt } \wTheta_0 = \Ov{\vr} \left( G +  |x_h|^2 \right) + {\rm const}.
	\end{equation}
	Assuming, for instance,  
	\[
	\avintO{ \wtilde{\mathcal{R}}_0 } = 0, 
	\]
	the density distribution $\wtilde{\mathcal{R}}_0$ is uniquely determined by $\wTheta_0$.
		
\end{Remark}

\medskip
One observes that, weak solutions of the synchronized system are defined in a similar way to Definition \ref{DL1} with obvious modifications reflecting the change of the initial 
time and the extra nudging terms in the field equations. In particular, the weak formulation of the momentum equation reads 
	\begin{align}
	\int_{0}^{T^+} &\intO{ \left[ \tvr \tvu \cdot \partial_t \bfphi + \tvr \tvu \otimes \tvu : \Grad \bfphi - {\frac{1}{\sqrt{\ep}} (\bf{e}_3 \times \tvr \tvu) \cdot 
			\bfphi} + 
		\frac{1}{\ep^2} p(\tvr, \tvt) \Div \bfphi \right] } \dt \br &= \int_{0}^{T^+} \intO{ \left[ \mathbb{S}(\tvt, \Ds \tvu) : \Grad \bfphi - \frac{1}{\ep} \tvr \Grad G \cdot \bfphi - { \frac{1}{\ep} \tvr \Grad |{x}_h |^2 \cdot \bfphi }  \right] } \dt \br
	&+ \Lambda \int_0^{T^+} \mathds{1}_{t \in [0,T]}\intO{ \Big( I_\delta[\tvu] - I_{\delta}[\vu] \Big)  \cdot \bfphi    }\dt  \br
	&- 
	\intO{ \tvr_{0} \tvu_{0} \cdot \bfphi (0, \cdot) }
	\label{m17}
\end{align}	
for any $\bfphi \in C^1_c([0, T^+) \times \Ov{\Omega}; R^3)$ such that $\bfphi|_{\partial B(r)} = 0$, 
$\bfphi \cdot \vc{n}|_{x_3 = 0,1} = 0$. The entropy production rate $\wtilde{\sigma}$ satisfies 
	\begin{equation} \label{m18}
	\wtilde{\sigma} \geq \frac{1}{\tvt} \left( \ep^2 \mathbb{S}(\tvt, \Ds \tvu) : \Ds \tvu - \frac{\vc{q}(\tvt ,\Grad \tvt): \Grad \tvt}{\tvt}\right) 
	- \Lambda \Big( I_\delta [\tvt] - I_\delta [\vt] \Big) \mathds{1}_{t \in [0,T]}.
\end{equation}
Finally, the ballistic energy balance \eqref{m6} is replaced by
\begin{align}  
	- &\int_{0}^{T^+} \partial_t \psi	\intO{ \left[ \ep^2 \frac{1}{2} \tvr |\tvu|^2 + \tvr e(\tvr, \tvt) - \Theta \tvr s(\tvr, \tvt) \right] } \dt  \br 
	&+ \int_{0}^{T^+} \!\! \int_{\Ov{\Omega}} \psi \,  \Theta\  \D \wtilde{\sigma} (t,x) + 
	\Lambda \int_0^{T^+} \psi \mathds{1}_{t \in [0,T]} \intO{ \tvt  \Big(  I_\delta[\tvt] - I_\delta [\vt]           \Big)        } \dt \br
	&+ \Lambda \int_0^{T^+} \psi \mathds{1}_{t \in [0,T]} \intO{  \Big( I_\delta[\tvu] - I_\delta [\vu] \Big) 
		\cdot \tvu  } \dt 
	\br
	&\leq 
	\int_{0}^{T^+} \psi \intO{ \left[ \ep \tvr \tvu \cdot \Grad G + { \ep \tvr \tvu \cdot \Grad |{x}_h |^2 } \right] } \dt \br  &- 
	\int_{0}^{T^+} \psi \intO{ \left[ 
		\tvr s(\tvr, \tvt) \partial_t \Theta + \tvr s(\tvr, \tvt) \tvu \cdot \Grad \Theta  + \frac{\vc{q}(\tvt, \Grad \tvt)}{\tvt} \cdot \Grad \Theta \right] } \dt \br 
	&+ \psi(0) \intO{ \left[ \frac{1}{2} \ep^2 \tvr_{0} |\tvu_{0}|^2 + \tvr_{0} e(\tvr_{0}, \tvt_{0}) - \Theta (0, \cdot) \tvr_{0} s(\tvr_{0}, \tvt_{0}) \right] }
	\label{m19}
\end{align}
for any $\psi \in C^1_c [0, T^+)$, $\psi \geq 0$, and \emph{any} smooth extension $\Theta$ of the boundary temperature, 
\[
\Theta > 0 \ \mbox{in}\ [0,T^+] \times \Ov{\Omega},\ \Theta|_{\partial \Omega} = \Ov{\vt} + \ep \vtB.                          
\]

\subsection{Main results}

Having collected the necessary preliminary material, we are ready to state our main result. 
First, let us introduce the norm of the \emph{data}, 
\begin{align}
\| \data \| &:= \Ov{\vr} + \Ov{\vr}^{-1} + \Ov{\vt} + \Ov{\vt}^{-1} + \| \vtB \|_{W^{2, \infty}(\partial \Omega)} 
\br
&\quad + \| \widetilde{\mathcal{T}}_0 \|_{W^{2, \infty}(\Omega)} 
+ \| {\mathcal{T}}_0 \|_{W^{2, \infty}(\Omega)} 
\br
&\quad + \| \widetilde{\vu}_{0,h} \|_{W^{2,\infty}(B(r); R^2)} 
+ \| {\vu}_{0,h} \|_{W^{2,\infty}(B(r); R^2)}.
\nonumber
\end{align}

\begin{Theorem}[\bf Tracking property] \label{mT1}
	
Suppose the thermodynamic functions $p$, $e$, and $s$ satisfy the hypotheses \eqref{w9}--\eqref{w14a}, and the transport coefficients 
$\mu$, $\eta$, and $\kappa$ are continuously differentiable functions of the temperature satisfying \eqref{w16}.

Let $(\vr, \vt, \vu)$ be a weak solution of the observed NSF system \eqref{s2}--\eqref{s4} in 
$(T^-, T^+) \times \Omega$, with the boundary data \eqref{s9}--\eqref{s11}, and the initial data satisfying 
\eqref{s12}--\eqref{s14}. Let $(\tvr, \tvt ,\tvu)$ be a weak solution of the synchronized NSF system \eqref{m7}--\eqref{m9} 
in $(0, T^+) \times \Omega$, with the boundary data \eqref{c5}--\eqref{m11}, the initial data \eqref{m12}--\eqref{m14}, and 
the interpolation operators $I_\delta$ satisfying \eqref{m10}.

Finally, let $0 < T < T^+$ and $\omega > 0$ be given. 

Then there exist
\begin{equation} \label{mm1}
\Lambda = 
\Lambda  \Big( \| \data \|, T^-, T^+ , T, \omega \Big) > 0, \ \delta = \delta ( \| \data \|, T^-, T^+ , T, \omega  ) > 0,
\end{equation}
and 
\begin{equation} \label{mm1a} 
  \ep_0 = \ep_0 \Big( \| \data \|, T^-, T^+ , T , \omega \Big) > 0,\ d_0 = d_0 \Big( \| \data \|, T^-, T^+ , T , \omega \Big) > 0
\end{equation}
such that
\begin{equation} \label{mm2}
\sup_{t \in [T, T^+]} \left( \left\| \left( \frac{\vr - \tvr}{\ep} \right)(t, \cdot) \right\|_{L^1(\Omega)} + 
\left\| \left( \frac{\vt - \tvt}{\ep} \right)(t, \cdot) \right\|_{L^1(\Omega)} + 
\| (\vr \vu - \tvr \tvu)(t, \cdot) \|_{L^1(\Omega; R^3)} \right) < \omega
\end{equation}
whenever 
\begin{equation} \label{mm3}
0 < \ep < \ep_0,\ 0 < d < d_0.
\end{equation}

\end{Theorem}

Note carefully that the initial states of the observed and the synchronized solutions are arbitrary, in particular they may be completely different. In the real world applications, 
the observed initial state is {\it a priori} unknown, while the initial data for the synchronized system can be 
conveniently chosen for (easy) numerical implementation, 
cf. Remark \ref{cR1}. The essential piece of information is the proximity of $(\vr, \vt, \vu)$ and $(\tvr, \tvt, \tvu)$ on the 
time interval $(T, T^+)$, where they satisfy the same system of equations and the synchronized solution provides a computable \emph{prediction} of the 
time evolution of the system. Finally, note that neither the observed solution $(\vr, \vt, \vu)$ nor its synchronized counterpart $(\tvr, \tvt, \tvu)$ 
are (known to be) uniquely determined by the data.

\medskip
The rest of the paper is devoted to the proof of Theorem \ref{mT1}.

\section{Asymptotic limit of the observed NSF system as $\ep \to 0$}
\label{as}

We start by performing the singular limit $\ep \to 0$, $d \to 0$ in the observed system. As shown in \cite{FanFei2025}, 
the asymptotic limit solves the \emph{target system} 
\begin{align} 
\Divh \vuh &= 0, \br	
\Ov{\vr} \Big[ \partial_t \vuh + \Divh (\vuh \otimes \vuh) \Big] + \Gradh \Pi &= 
\mu(\Ov{\vt}) \Delta_h \vuh + \left<{\mathcal{R}} \right> \Gradh \Big( G + \frac{1}{2}|x_h|^2 \Big) \ \mbox{in}\ (T^-,T^+) \times B(r), \br 
\Ov{\vr} c_p(\Ov{\vr}, \Ov{\vt} ) \Big[ \partial_t \mathcal{T} + \vuh \cdot \Gradh \mathcal{T} \Big] &- 
\Ov{\vr} \Ov{\vt} \alpha (\Ov{\vr}, \Ov{\vt} ) \vuh \cdot \Gradh \Big( G + \frac{1}{2}| x_h|^2 \Big) \br &= 
\kappa (\Ov{\vt}) \Del \mathcal{T} + \Ov{\vt}  \alpha (\Ov{\vr}, \Ov{\vt} )  \frac{\partial p  (\Ov{\vr}, \Ov{\vt} ) } 
{\partial \vt} \partial_t \avintO{ \mathcal{T} } \ \mbox{in}\ (T^-, T^+) \times \Omega, 
\label{a2}
\end{align}
where $\left<{\mathcal{R}} \right>:= \int_0^1 {\mathcal{R}} \dx_3$ and supplemented with the Boussinesq relation 
\begin{equation} 
	\frac{\partial p(\Ov{\vr}, \Ov{\vt})}{\partial \vr} \Grad {\mathcal{R}} + \frac{\partial p(\Ov{\vr}, \Ov{\vt})}{\partial \vt} \Grad \mathcal{T} =
\Ov{\vr} \Big( \Grad G + \frac{1}{2} \Grad |{x}_h|^2 \Big),\ \intO{ {\mathcal{R}} } = 0,
\label{aa2}	
\end{equation}	
with the boundary conditions
\begin{equation} \label{a3}
\vuh|_{\partial B(r)} = 0,\ 
\mathcal{T}|_{\partial \Omega} = \vtB, 
\end{equation}	
and the initial data 
\begin{equation} \label{a4}
\vuh(T^-, \cdot) = {\vu}_{0,h},\ \mathcal{T}(T^-,\cdot) = {\mathcal{T}}_0.
\end{equation}
Here, $c_p(\Ov{\vr}, \Ov{\vt})$ and $\alpha(\Ov{\vr}, \Ov{\vt})$ stand for the specific heat at constant pressure and the thermal expansion coefficient, respectively, evaluated at the reference state $(\Ov{\vr}, \Ov{\vt})$, see  \cite{FanFei2025} for details. A remarkable feature of system \eqref{a2} is 
that the fluid motion is purely horizontal, governed by the 2-D incompressible Navier--Stokes equations. In particular, problem \eqref{a2}--\eqref{a4} 
is globally well posed as soon as the data are smooth enough, see \cite{FeGwSG25}. 

\medskip
The main result of \cite[Theorem 3.1]{FanFei2025} asserts that for any $\omega^1 > 0$, there exist 
$\ep^1_0 > 0$, $d^1_0 > 0$ depending on $\| \data \|$, notably the norms of $\vrB$, $\mathcal{T}_0$, $\vu_{0,h}$, and the length of the time interval 
$T^+ - T^-$, such that 
\begin{align}
\sup_{t \in [T^-, T^+]} \left( \left\| \left( \frac{ \vr - \Ov{\vr} }{\ep} - {\mathcal{R}} \right)(t, \cdot) \right\|_{L^1(\Omega)}
+ \left\| \left( \frac{ \vt - \Ov{\vt} }{\ep} - {\mathcal{T}} \right)(t, \cdot) \right\|_{L^1(\Omega)} \right) &< \omega^1, 
\br 
\sup_{t \in [T^-, T^+]} \left( \left\| \left( \vr \vu - \Ov{\vr} \vuh \right)(t, \cdot) \right\|_{L^1(\Omega; R^3)}  \right) &< \omega^1, 
\label{a5} \\
\left\| \left( \frac{ \vt - \Ov{\vt} }{\ep} - {\mathcal{T}} \right) \right\|_{L^2(T^-,T^+; W^{1,2}(\Omega))} 
+ \left\| \vu - \vuh \right\|_{L^2(T^-,T^+; W^{1,2}(\Omega; R^3))} &< \omega^1 
\label{a6}  
\end{align}
whenever $0 < \ep < \ep^1_0$, $0 < d < d^1_0$ for \emph{any} weak solution $(\vr, \vt, \vu)$ of the observed NSF system with the initial data \eqref{s12}--\eqref{s14}.

\begin{Remark} \label{asR1}
As a matter of fact, the main result of \cite{FanFei2025} is stated in a slightly different setting replacing the time derivative of the temperature 
average in \eqref{a2} by a non--local boundary condition.
\end{Remark}	

\section{A data assimilation method for the target system}
\label{das}

Our next goal is to adapt the assimilation data method of Azouani, Olson, and Titi \cite{AzOlTi} to the target system \eqref{aa2}. 
Specifically, we consider a synchronized problem 
\begin{align} 
	\Divh \wvuh &= 0, \br	
	\Ov{\vr} \Big[ \partial_t \wvuh + \Divh (\wvuh \otimes \wvuh) \Big] + \Gradh \wtilde{\Pi} &= 
	\mu(\Ov{\vt}) \Delta_h \wvuh + \left<\wtilde{\mathcal{R}} \right> \Gradh \Big( G + \frac{1}{2}|\vc x_h|^2 \Big) 
	- \Lambda I_{\delta} [ \wvuh - \vuh] \mathds{1}_{t \in [0,T]} , \label{A4} \\ 
	\mbox{in}\ (0,T^+) \times B(r), \br  \Ov{\vr} c_p(\Ov{\vr}, \Ov{\vt} ) \Big[ \partial_t \wTheta + \wvuh \cdot \Gradh \wTheta \Big] &- 
	\Ov{\vr} \Ov{\vt} \alpha (\Ov{\vr}, \Ov{\vt} ) \wvuh \cdot \Gradh \Big( G + \frac{1}{2}|\vc x_h|^2 \Big) \br &= 
	\kappa (\Ov{\vt}) \Del \wTheta + \Ov{\vt}  \alpha (\Ov{\vr}, \Ov{\vt} )  \frac{\partial p  (\Ov{\vr}, \Ov{\vt} ) } 
	{\partial \vt} \partial_t \avintO{ \wTheta } - \Ov{\vt} \Lambda I_\delta \left[ \wTheta - \mathcal{T} \right]\mathds{1}_{t \in [0,T]} \label{A5}\\
	\mbox{in}\ (0, T^+) \times \Omega, \br 
	\frac{\partial p(\Ov{\vr}, \Ov{\vt})}{\partial \vr} \Grad \wtilde{\mathcal{R}} + \frac{\partial p(\Ov{\vr}, \Ov{\vt})}{\partial \vt} \Grad \wTheta &=
	\Ov{\vr} \Big( \Grad G + \frac{1}{2} \Grad |{x}_h|^2 \Big),\ \intO{ \wtilde{\mathcal{R}} } = 0,
	\label{A6}	
\end{align}	
with the boundary conditions
\begin{equation} \label{A7}
	\wvuh|_{\partial B(r)} = 0,\ 
	\wTheta|_{\partial \Omega} = \vtB, 
\end{equation}	
and the initial data 
\begin{equation} \label{A8}
	\wvuh(0, \cdot) = {\wtilde{\vu}}_{0,h},\ \wTheta(0,\cdot) = \wtilde{\mathcal{T}}_0.
\end{equation}

The nudging operators in \eqref{A4}, \eqref{A5} force the solution $(\wvuh, \wTheta)$ to remain on the proximity of the solution $(\vuh, \mathcal{T})$ 
of the target system \eqref{a2}--\eqref{a4} in the time interval $[T, T^+]$.

\subsection{Global existence for the target problem}

We claim that both the target system \eqref{a2}--\eqref{a4} and its synchronized counterpart \eqref{A4}--\eqref{A8} are well posed in the standard 
$L^p-$framework. 

\begin{Proposition} \label{dasP1} {\bf (Well posedness of the target system)}

\begin{itemize}

\item

The target system \eqref{a2}--\eqref{a4} admits a strong solution $(\vuh, \mathcal{T})$ in $(T^-, T^+) \times (B(r) \times \Omega)$ 
for any initial data $(\vu_{0,h}, \mathcal{T}_0)$ specified in \eqref{s13}--\eqref{s14}, unique in the class
\begin{align}
\vuh &\in L^q(T^-, T^+; W^{2,q} (B(r); R^2)),\ \partial_t \vuh \in L^q(T^-, T^+; L^q (B(r); R^2)) \ \mbox{for any}\ 1 \leq q < \infty, \br
\mathcal{T} &\in L^q(T^-, T^+; W^{2,q} (\Omega)),\ \partial_t \mathcal{T} \in L^q(T^-, T^+; L^q (\Omega)) \ \mbox{for any}\ 1 \leq q < \infty.
\label{das1}
\end{align}  	
\item The synchronized target system \eqref{A4}--\eqref{A8} admits a strong solution $(\wvuh, \wTheta)$ in $(0, T^+) \times (B(r) \times \Omega)$
for any initial data $(\widetilde{u}_{0,h}, \wTheta_0)$ specified in \eqref{m12}--\eqref{m14}, unique in the class
\begin{align}
	\wvuh &\in L^q(0, T^+; W^{2,q} (B(r); R^2)),\ \partial_t \wvuh \in L^q(0, T^+; L^q (B(r); R^2)) \ \mbox{for any}\ 1 \leq q < \infty, \br
	\wTheta &\in L^q(0, T^+; W^{2,q} (\Omega)),\ \partial_t \wTheta \in L^q(0, T^+; L^q (\Omega)) \ \mbox{for any}\ 1 \leq q < \infty.
	\label{das2}
\end{align}
for any fixed $\Lambda > 0$, $\delta > 0$.	
\end{itemize}

\end{Proposition}

The result for the target system was proved in \cite[Theorem 2.1]{FeGwSG25}. 
As the interpolation operators range in $L^\infty$, 
the proof can be adapted to the synchronized system in a direct manner. Note however that certain norms for the synchronized solution may depend on $\Lambda$ and $\delta$.

\subsection{Tracking property for the target system}

Our next objective is to show that the solution of the target problem and that of the synchronized problem stay close on the 
time interval $[T, T^+]$ as soon as $\Lambda > 0$ and $\delta > 0$ are conveniently fixed. To this end, we report the integral identities 
proved in \cite[formulae (4.13), (4.14)]{FeGwSG25}:
\begin{align} 
	\frac{1}{2} \frac{\D}{\dt} \| \vuh - \wvuh \|^2_{L^2(B(r); R^2)}  &+ \frac{\mu (\Ov{\vt})}{\Ov{\vr}} \| \nabla_h( \vuh - \wvuh) \|^2_{L^2(B(r); R^4)}
	+ \Lambda \intOh{ I_\delta[ \vuh - \wvuh ] \cdot (\vuh - \wvuh) } \mathds{1}_{t \in [0,T]} \br &= \frac{1}{2 \Ov{\vr}} \intOh{ \left< \wTheta - \mathcal{T} \right > \Gradh |x_h|^2 \cdot (\vuh - \wvuh) } \br
	&+ \intOh{ \Divh( \wvuh \otimes \wvuh - \vuh \otimes \vuh) \cdot (\vuh - \wvuh) },
	\label{das3}
\end{align}	
and 
\begin{align} 
	\frac{1}{2} &\frac{\D }{\dt} \left[ \intO{ (\mathcal{T} - \wTheta)^2 } - 
	\frac{\beta (\Ov{\vr}, \Ov{\vt})}{|\Omega|}  \left( \intO{ (\mathcal{T} - \wTheta) } \right)^2  \right]+ K(\Ov{\vr}, \Ov{\vt}) \intO{ |\Grad (\mathcal{T} - \wTheta) |^2 } \br
	&= \gamma (\Ov{\vr}, \Ov{\vt}) \intO{ \mathcal{T} (\vuh - \wvuh) \cdot \Grad (\mathcal{T} - \wTheta)  } - 
	\Lambda \intO{ I_\delta \Big[\mathcal{T} - \wTheta \Big] \Big( \mathcal{T} - \wTheta \Big) } \mathds{1}_{t \in [0,T]} , 
	\label{das4}
\end{align}	
where $0 < \beta (\Ov{\vr}, \Ov{\vt}) < 1$, $K (\Ov{\vr}, \Ov{\vt}) > 0$, and $\gamma = \gamma (\Ov{\vr}, \Ov{\vt})$ are constants depending solely on the 
reference state $(\Ov{\vr}, \Ov{\vt})$. 

Furthermore, thanks to solenoidality of both $\vuh$ and $\wvuh$, the last integral in \eqref{das3} can be written as
\[
\intOh{ \Divh( \wvuh \otimes \wvuh - \vuh \otimes \vuh) \cdot (\vuh - \wvuh) } 
= - \intOh{ (\vuh - \wvuh) \cdot \Grad \vuh \cdot (\vuh - \wvuh) }
\]
where, by Ladyzhenskaya interpolation inequality ,
\begin{align} 
&\left| \intOh{ (\vuh - \wvuh) \cdot \Grad \vuh \cdot (\vuh - \wvuh) } \right| 
\leq \| \nabla_h \vuh \|_{L^2(B(r); R^4)} \| \vuh - \wvuh \|^2_{L^4(\Omega_h; R^2)} \br &\leq C_L  \| \nabla_h \vuh \|_{L^2(B(r); R^4)} 
 \| \vuh - \wvuh \|_{L^2(B(r); R^2)} \| \nabla_h (\vuh - \wvuh) \|_{L^2(B(r); R^4)} .
\label{das5}
\end{align}

Similarly, 
\begin{equation} \label{das6}
\left| \intO{ \mathcal{T} (\vuh - \wvuh) \cdot \Grad (\mathcal{T} - \wTheta)  } \right| 
\leq \| \mathcal{T} \|_{L^\infty(\Omega)} \| \vuh - \wvuh \|_{L^4(\Omega_h; R^2)} 
\| \Grad (\mathcal{T} - \wTheta) \|_{L^2(\Omega; R^3)}.
\end{equation}
It is important that both $\| \nabla_h \vuh \|_{L^2(B(r); R^4)}$ in \eqref{das5} and $\| \mathcal{T} \|_{L^\infty(\Omega)}$ in \eqref{das6}
depend only on the observed solution, meaning on $\| \data \|$. In particular, they are independent of $\Lambda$ and $\delta$.  

Using the properties of the interpolation operators we obtain 
\begin{align}
\Lambda \intOh{ I_\delta [\vuh - \wvuh] \cdot (\vuh - \wvuh) } &= 
\Lambda \| \vuh - \wvuh \|^2_{L^2(B(r); R^2)}\br &+ \Lambda \intOh{ [I_\delta - \mathbb{I}] [\vuh - \wvuh] \cdot (\vuh - \wvuh) }, \br
\Lambda \intO{ I_\delta \Big[ \mathcal{T} - \wTheta \Big] \Big( \mathcal{T} - \wTheta \Big) } &= 
\Lambda \| \mathcal{T} - \wTheta \|_{L^2(\Omega)}^2 \br 
&+ \Lambda \intOh{ [I_\delta - \mathbb{I}] [\mathcal{T} - \wTheta] \cdot (\mathcal{T} - \wTheta) },
\label{das7}
\end{align}
where, by virtue of \eqref{m10}, 
\begin{align}
\left|	\intOh{ [I_\delta - \mathbb{I}] [\vuh - \wvuh] \cdot (\vuh - \wvuh) } \right| 
&\leq \delta \| \vuh - \wvuh \|_{W^{1,2}_0(B(r); R^2)}\| \vuh - \wvuh \|_{L^2(B(r); R^2)}, \br
\left|	\intOh{ [I_\delta - \mathbb{I}] [\mathcal{T} - \wTheta] \cdot (\mathcal{T} - \wTheta) } \right|
&\leq \delta \| \mathcal{T} - \wTheta \|_{W^{1,2}_0(\Omega)}\| \mathcal{T} - \wTheta \|_{L^2(\Omega)}.
\label{das8}	
\end{align}

Next, we recall Poincar\' e inequality
\begin{align}
\| \vuh - \wvuh \|_{W^{1,2}_0(B(r); R^2)} &\leq C_P \| \nabla_h (\vuh - \wvuh) \|_{L^{2}(B(r); R^4)}, \br
\| \mathcal{T} - \wTheta \|_{W^{1,2}_0(\Omega)} &\leq C_P \| \Grad (\mathcal{T} - \wTheta) \|_{L^{2}(\Omega; R^3)}.
\label{das9}	
\end{align}

Summing the integral identities \eqref{das3}, \eqref{das4} with the estimates \eqref{das5}--\eqref{das9}, we 
deduce there exist $\Lambda_0 = \Lambda_0 (\| \data \|) > 0$ and $\delta_0 = \delta_0(\Lambda_0) > 0$ such that 
\begin{align}
	\frac{\D}{\dt} &\left( \| \vuh - \wvuh \|^2_{L^2(B(r); R^2)} + \left[ \intO{ (\mathcal{T} - \wTheta)^2 } - 
	\frac{\beta (\Ov{\vr}, \Ov{\vt})}{|\Omega|}  \left( \intO{ (\mathcal{T} - \wTheta) } \right)^2  \right]  \right) \br
	&\leq - \frac{\Lambda}{2} \left(  \| \vuh - \wvuh \|^2_{L^2(B(r); R^2)} + 
	\| \mathcal{T} - \wTheta \|^2_{L^2(\Omega)} \right)
\label{das10}	
\end{align}
for all $t \in [0,T]$ whenever $\Lambda > \Lambda_0$, $0 < \delta = \delta(\Lambda) \leq \delta_0$. By the same token, 
there exists $K = K(\|\data\|)$ such that 
\begin{align}
	\frac{\D}{\dt} &\left( \| \vuh - \wvuh \|^2_{L^2(B(r); R^2)} + \left[ \intO{ (\mathcal{T} - \wTheta)^2 } - 
	\frac{\beta (\Ov{\vr}, \Ov{\vt})}{|\Omega|}  \left( \intO{ (\mathcal{T} - \wTheta) } \right)^2  \right]  \right) \br
	&\leq K(\|\data\|) \left(  \| \vuh - \wvuh \|^2_{L^2(B(r); R^2)} + 
	\| \mathcal{T} - \wTheta \|^2_{L^2(\Omega)} \right)
	\label{das11}	
\end{align}
for all $t \in (T, T^+)$. Finally, by Jensen's inequality
\begin{equation} \label{das12}
\Big(1 - \beta(\Ov{\vr}, \Ov{\vt}) \Big)	\| \mathcal{T} - \wTheta \|_{L^2(\Omega)}^2 \leq
\left[ \intO{ (\mathcal{T} - \wTheta)^2 } - 
\frac{\beta (\Ov{\vr}, \Ov{\vt})}{|\Omega|}  \left( \intO{ (\mathcal{T} - \wTheta) } \right)^2  \right] \leq 
\| \mathcal{T} - \wTheta \|_{L^2(\Omega)}^2. 
\end{equation}
Combining \eqref{das10}--\eqref{das12} we may infer that for any $\omega^2 > 0$ there exist 
\[
\Lambda = \Lambda (\| \data \|, T, T^+, \omega^2) > 0,\ 
\delta = \delta (\| \data \|, T, T^+, \omega^2) > 0,  
\]
such that
\begin{equation} \label{das13}
\sup_{t \in [T, T^+ ] } \left( \| (\mathcal{T} - \wTheta)(t, \cdot) \|_{L^2(\Omega)} + \| (\mathcal{R} - \wtilde{\mathcal{R}})(t, \cdot) \|_{L^2(\Omega)} 
+ \| (\vuh - \wvuh)(t, \cdot) \|^2_{L^2(B(r); R^2)}
\right) < \omega_2. 
\end{equation}	

\section{Asymptotic limit of the synchronized NSF system as $\ep \to 0$}
\label{F}

With the parameters $\Lambda > 0$ and $\delta > 0$ fixed, we perform the asymptotic limit in the 
synchronized NSF system \eqref{m7}--\eqref{m14}. Specifically, we consider a family of 
synchronized solutions $(\tvre, \tvte, \tvue)_{\ep >0 }$, with the associate observed solutions $(\vre, \vte, \vue)_{\ep > 0}$ 
appearing in the ``nudging'' integrals, and perform the asymptotic limit $\ep \to 0$. By analogy with 
the results of \cite{FanFei2025}, a natural candidate for the asymptotic limit is the unique solution $(\wtilde{R}, \wTheta, 
\wvuh)$
of the synchronized target system \eqref{A4}--\eqref{A8}.

Let us start by introducing the relative energy functional 
\begin{align} 
E_\ep &\left( \tvr, \tvt, \tvu \Big| r, \Theta, \vU \right) = \frac{1}{2} \tvr |\tvu - \vU |^2 \br 
&+\frac{1}{\ep^2} \left[ \tvr e(\tvr, \tvt) - \Theta \Big( \tvr s(\tvr , \tvt) - r s(r, \Theta) \Big) - \left( e(r, \Theta) - \Theta s(r, \Theta) + \frac{p(r, \Theta)}{\Theta} \right)(\vr - r) - r e(r, \Theta)    \right].	
\label{A1}
\end{align}
As the total energy $E_\ep$ is a (strictly) convex function of the density, total entropy, and momentum, the relative energy is nothing other than the 
associated Bregmann distance, see \cite{FeiNovOpen} for details.

Let $(\tvr, \tvt, \tvu)$ be a weak solution of the synchronized NSF system \eqref{m7}--\eqref{m14}. Note that the \emph{existence} 
can be shown in the same way as in \cite[Chapter 12]{FeiNovOpen} as the additional ``damping'' terms do not bring any essential 
difficulty (thanks to our choice of $I_\delta$). 
First, we recall the relative energy inequality (see \cite{FanFei2025}) augmented by the nudging terms:
 \begin{align}
	&\left[ \intO{ E_\ep \left(\tvr, \tvt, \tvu \Big| r, \Theta, \vU \right) } \right]_{t = 0}^{t = \tau} \br 
	&+ \int_0^\tau  \int_{{\Omega}} \frac{\Theta}{\tvt} \left(  \mathbb{S}(\tvt, \Ds \tvu) : \Ds \tvu - \frac{1}{\ep^2} \frac{\vc{q}(\tvt ,\Grad \tvt)\cdot \Grad \tvt}{\tvt}\right)   \dx \dt \br 
	&\leq - \frac{1}{\ep^2} \int_0^\tau \intO{ \left( \tvr (s(\tvr, \tvt) - s(r, \Theta)) \partial_t \Theta + \tvr (s(\tvr,\tvt) - s(r, \Theta)) \tvu \cdot \Grad \Theta \right) } \dt \br 
	&+ \frac{1}{\ep^2} \int_0^\tau \intO{ 
		\frac{\kappa (\tvt) \Grad \tvt}{\tvt}  \cdot \Grad \Theta } \dt \br 
	&- \int_0^\tau \intO{ \Big[ \tvr (\tvu - \vU) \otimes (\tvu - \vU) + \frac{1}{\ep^2} p(\tvr, \tvt) \mathbb{I} - \mathbb{S}(\tvt, \Ds \tvu) \Big] : \Ds \vU } \dt   \br
	&+ \int_0^\tau \intO{ \tvr \left[ \frac{1}{\ep} \Grad G + \frac{1}{2\ep} \Grad |{x}_h |^2 
		- \frac{2}{\sqrt{\ep}} (\vc{e}_3 \times \tvu)  - \partial_t \vU - (\vU \cdot \Grad) \vU  \right] \cdot (\tvu - \vU) } \dt  \br 
	&+ \frac{1}{\ep^2} \int_0^\tau \intO{ \left[ \left( 1 - \frac{\tvr}{r} \right) \partial_t p(r, \Theta) - \frac{\tvr}{r} \tvu \cdot \Grad p(r, \Theta) \right] } \dt \br 
	& - \Lambda \int_0^\tau \mathds{1}_{t \in [0,T] } \intO{ I_{\delta}[ \tvu - \vu ] \cdot (\tvu - \vU) } \dt - 
	\Lambda \int_0^\tau \mathds{1}_{t \in [0,T] }  \intO{ \left( \frac{\tvt - \Theta}{\ep} \right) I_\delta\left[ \frac{\tvt - \vt}{\ep} \right] } \dt,  
	\label{A2}
\end{align}
for a.a. $0 \leq \tau < T^+$ and any trio of continuously differentiable functions $(r, \Theta, \vU)$ satisfying
\begin{equation} \label{A3}
	r > 0,\quad \Theta > 0,\quad \Theta|_{\partial \Omega} = \Ov{\vt} + \ep \vtB, \quad \vU|_{\partial B(r) \times (0,1)} = 0,\quad \vU \cdot \vc{n}|_{x_3 = 0,1} = 0.
\end{equation}

\subsection{Ansatz in the relative energy inequality}

Similarly to \cite{FanFei2025}, we consider the relative energy inequality \eqref{A2} for 
\begin{equation} \label{ansatz}
(\tvr, \tvt, \tvu) = (\tvre, \tvte, \tvue), \ (\vt, \vu) = (\vte, \vue), \  r = \Ov{\vr} + \ep \wtilde{\mathcal{R}},\ 
 \Theta = \Ov{\vt} + \ep \wTheta,\ \vU = \wvuh. 
\end{equation}
Thanks to our choice of the initial data \eqref{A8} for the synchronized problems \eqref{m15} and \eqref{A8}, 
we have 
\begin{align} \label{F1}
&\intO{ E_\ep \left( \tvre, \tvte, \tvue \Big| \Ov{\vr} + \ep \wtilde{\mathcal{R}} , \Ov{\vt} + \ep \wtilde{\mathcal{T}}, 
\wvuh \right)(0, \cdot) } \br &\quad = 	\intO{ E_\ep \left( \tvre(0, \cdot), \tvte (0, \cdot), \tvue (0, \cdot) \Big| \Ov{\vr} + \ep \wtilde{\mathcal{R}}_0 , \Ov{\vt} + \ep \wtilde{\mathcal{T}}_0, 
\wtilde{\vu}_{0,h}  \right) } \to 0 \ \mbox{as}\ \ep \to 0. 
\end{align}

Our final goal is to show 
\begin{equation} \label{F2}
\intO{ E_\ep \left( \tvre, \tvte, \tvue \Big| \Ov{\vr} + \ep \wtilde{\mathcal{R}} , \Ov{\vt} + \ep \wtilde{\mathcal{T}}, 
	\wvuh \right)(\tau, \cdot) } \to 0 \ \mbox{as}\ \ep \to 0, d \to 0 \ \mbox{uniformly for}\ \tau \in [0,T].
\end{equation}	
Relation \eqref{F2}, once proved, will yield the following conclusion. For any $\omega^3 > 0$, there exist
$\ep^2_0 = \ep^2_0 (\| \data \|,  T^+, T, \omega^3) > 0$, 
$d^2_0 = d^2_0 (\| \data \|,  T^+, T, \omega^3) > 0$ so that 
\begin{align}
	\sup_{t \in [0, T^+]} \left( \left\| \left( \frac{ \tvr - \Ov{\vr} }{\ep} - \wtilde{\mathcal{R}} \right)(t, \cdot) \right\|_{L^1(\Omega)}
	+ \left\| \left( \frac{ \tvt - \Ov{\vt} }{\ep} - \wtilde{\mathcal{T}} \right)(t, \cdot) \right\|_{L^1(\Omega)} \right) &< \omega^3, 
	\br 
	\sup_{t \in [0, T^+]} \left( \left\| \left( \tvr \tvu - \Ov{\vr} \wvuh \right)(t, \cdot) \right\|_{L^1(\Omega; R^3)}  \right) &< \omega^3, 
	\label{F3} \\
	\left\| \left( \frac{ \tvt - \Ov{\vt} }{\ep} - \wtilde{\mathcal{T}} \right) \right\|_{L^2(0,T^+; W^{1,2}(\Omega))} 
	+ \left\| \tvu - \wvuh \right\|_{L^2(0,T^+; W^{1,2}(\Omega; R^3))} &< \omega^3 
	\label{F4}  
\end{align}
whenever $0 < \ep < \ep^2_0$, $0 < d < d^2_0$ for \emph{any} weak solution $(\tvr, \tvt, \tvu)$ of the synchronized NSF system with the initial data \eqref{m12}--\eqref{m14}. 

Now, combining the bounds \eqref{a5}, \eqref{a6} with \eqref{das13} and 
\eqref{F3}, \eqref{F4}, we complete the proof of Theorem \ref{mT1}. Indeed going back to \eqref{mm2} 
we obtain
\begin{align}
	&\left\| \left( \frac{\vr - \tvr}{\ep} \right)(t, \cdot) \right\|_{L^1(\Omega)} + 
	\left\| \left( \frac{\vt - \tvt}{\ep} \right)(t, \cdot) \right\|_{L^1(\Omega)} + 
	\| (\vr \vu - \tvr \tvu)(t, \cdot) \|_{L^1(\Omega; R^3)} \br
&\leq 	\left\| \left( \frac{\vr - \Ov{\vr}}{\ep} - \mathcal{R} \right)(t, \cdot) \right\|_{L^1(\Omega)} + 
\left\| \left( \frac{\vt - \Ov{\vt}}{\ep} - \mathcal{T} \right)(t, \cdot) \right\|_{L^1(\Omega)} + 
\| (\vr \vu - \Ov{\vr} \vuh)(t, \cdot) \|_{L^1(\Omega; R^3)} \br 
&+ \left\| (\mathcal{R} - \wtilde{\mathcal{R}} )(t, \cdot) \right\|_{L^1(\Omega)} + 
\left\| (\mathcal{T} - \wtilde{\mathcal{T}} )(t, \cdot) \right\|_{L^1(\Omega)} + 
\Ov{\vr} \| (\vuh - \wvuh)(t, \cdot) \|_{L^1(B(r); R^2)} \br 
&+ \left\| \left( \frac{\tvr - \Ov{\vr}}{\ep} - \wtilde{\mathcal{R}} \right)(t, \cdot) \right\|_{L^1(\Omega)} + 
\left\| \left( \frac{\tvt - {\Ov{\vt}}}{\ep} - \wtilde{\mathcal{T}} \right)(t, \cdot) \right\|_{L^1(\Omega)} + 
\| (\tvr \tvu - \Ov{\vr} \wvuh)(t, \cdot) \|_{L^1(\Omega; R^3)}.
\label{F5}
\end{align}
Thus given $\omega > 0$ we may choose $\omega^i = \omega/3$, $i = 1,2,3$, and 
$\ep_0 = \min\{ \ep^1_0, \ep^2_0 \}$, $d_0 =  \min\{ d^1_0, d^2_0 \}$ to obtain the desired conclusion \eqref{mm2}. The rest of this section is therefore 
devoted to the proof of \eqref{F2}.

\subsection{Additional terms in the relative energy inequality}

The relative energy inequality \eqref{A2} with the ansatz \eqref{ansatz} reads
 \begin{align}
	&\intO{ E_\ep \left(\tvre, \tvte, \tvue \Big| \Ov{\vr} + \ep \wtilde{\mathcal{R}} , \Ov{\vt} + \ep \wtilde{\mathcal{T}}, \wvuh \right) (\tau, \cdot) }  \br 
	&+ \int_0^\tau  \int_{{\Omega}} \frac{\Ov{\vt} + \ep \wtilde{\mathcal{T}} }{\tvte} \left(  \mathbb{S}(\tvte, \Ds \tvue) : \Ds \tvue + \frac{1}{\ep^2} \frac{\kappa (\tvte) \Grad \tvte \cdot \Grad \tvte}{\tvte}\right)   \dx \dt \br 
	&\leq - \frac{1}{\ep} \int_0^\tau \intO{ \tvre \Big( s(\tvre, \tvte) - s(\Ov{\vr} + \ep \wtilde{\mathcal{R}}, \Ov{\vt} + \ep \wtilde{\mathcal{T}}) \Big) \partial_t \wtilde{\mathcal{T}} } \dt \br 
	&\quad - \frac{1}{\ep} \int_0^\tau \intO{ \tvre \Big( s(\tvre, \tvte) - s(\Ov{\vr} + \ep \wtilde{\mathcal{R}}, \Ov{\vt} + \ep \wtilde{\mathcal{T}}) \Big) \tvue \cdot \Grad \wtilde{\mathcal{T}}  } \dt
	\br 
	&\quad + \frac{1}{\ep} \int_0^\tau \intO{ 
		\frac{\kappa (\tvte) \Grad \tvte}{\tvte}  \cdot \Grad \wtilde{\mathcal{T}} } \dt \br 
	&\quad - \int_0^\tau \intO{ \Big[ \tvre (\tvue - \wvuh) \otimes (\vue - \wvuh)  - \mathbb{S}(\tvte, \Ds \tvue) \Big] : \Ds \wvuh } \dt   \br
	&\quad + \int_0^\tau \intO{ \tvre \left[ \frac{1}{\ep} \Grad G + \frac{1}{2\ep} \Grad |{x}_h |^2 
		- \partial_t \wvuh - (\wvuh \cdot \Grad) \wvuh  \right] \cdot (\tvue - \wvuh) } \dt  \br 
	&\quad + \frac{1}{\ep^2} \int_0^\tau \intO{ \left[ \left( 1 - \frac{\tvre}{\Ov{\vr} + \ep \wtilde{\mathcal{R}}} \right) 
		\partial_t p(\Ov{\vr} + \ep \wtilde{\mathcal{R}}, \Ov{\vt} + \ep \wtilde{\mathcal{T}}) - \frac{\tvre}{\Ov{\vr} + \ep \wtilde{\mathcal{R}}} \tvue \cdot \Grad p(\Ov{\vr} + \ep \wtilde{\mathcal{R}}, \Ov{\vt} + \ep \widetilde{\mathcal{T}}) \right] } \dt\br 
	&\quad - \Lambda \int_0^\tau \mathds{1}_{t \in [0,T] } \intO{ I_{\delta}[ \tvue - \vue ] \cdot (\tvue - \wvuh) } \dt\br &\quad - 
	\Lambda \int_0^\tau \mathds{1}_{t \in [0,T] }  \intO{ \left( \frac{\tvte - \Ov{\vt}}{\ep} - 
		\wTheta \right) I_\delta\left[ \frac{\tvte - \vte}{\ep} \right] } \dt  + \mathcal{O}(\ep, \tau) ,
	\label{F6}
\end{align}
where the symbol $\mathcal{O}(\ep, \tau)$ denotes a generic function satisfying
\begin{equation} \label{def:omega_ep}
	\sup_{\tau\in[0,T^+]}\mathcal{O}(\ep, \tau) \to 0 \ \mbox{as}\ \ep \to 0,
\end{equation}
cf. \cite[Section 5.2.1]{FanFei2025}.

To show \eqref{F2}, we follow step by step the arguments of \cite{FanFei2025} estimating the integrals on the 
right--hand side by the left--hand side and finally applying a Gr\" onwall type argument. 
The main novelty with respect to
\cite{FanFei2025} lies in 
the presence of ``nudging terms'', which are carefully handled in the following section.

\subsubsection{Velocity perturbation}
\label{vp}

There are two extra terms on the right--hand side of \eqref{F6} related to the nudging integrals acting on the velocity fields. 
Besides 
\[
- \Lambda \int_0^\tau \mathds{1}_{t \in [0,T] } \intO{ I_{\delta}[ \tvue - \vue ] \cdot (\tvue - \wvuh) } \dt, 
\]
there is 
\begin{align}
- &\int_0^\tau \intO{ \tvre \left[  
\partial_t \wvuh + (\wvuh \cdot \Grad) \wvuh  \right] \cdot (\tvue - \wvuh) } \dt \br &= 
- \ep \int_0^\tau \intO{ \frac{\tvre - \Ov{\vr}}{\ep} \left[  
\partial_t \wvuh + (\wvuh \cdot \Grad) \wvuh  \right] \cdot (\tvue - \wvuh) } \dt \br 
&\quad -\int_0^\tau \intO{ \Ov{\vr} \left[  
\partial_t \wvuh + (\wvuh \cdot \Grad) \wvuh  \right] \cdot (\tvue - \wvuh) } \dt, 
\nonumber
\end{align}
where 
\[
-\int_0^\tau \intO{ \Ov{\vr} \left[  
\partial_t \wvuh + (\wvuh \cdot \Grad) \wvuh  \right] \cdot (\tvue - \wvuh) } \dt = \dots + \Lambda 
\int_0^\tau \mathds{1}_{t \in [0,T]} \intO{ I_\delta[ \wvuh - \vuh] \cdot (\tvue - \wvuh) } \dt.
\]
Thus 
the velocity perturbation yields a term on the right--hand side of the relative energy inequality in the form 
\begin{align}
- &\Lambda \int_0^\tau \mathds{1}_{t \in [0,T]} \intO{ I_{\delta}[ \tvue - \vue ] \cdot (\tvue - \wvuh) } \dt + 
\Lambda \int_0^\tau \mathds{1}_{t \in [0,T]}\intO{ I_{\delta}[ \wvuh - \vuh] \cdot  (\tvue - \wvuh) } \dt \br
&= - \Lambda \int_0^\tau \mathds{1}_{t \in [0,T]} \intO{ (\tvue - \wvuh) \cdot I_{\delta}[ \tvue - \wvuh ] }\dt + 
\Lambda \int_0^\tau \mathds{1}_{t \in [0,T]} \intO{ (\tvue - \wvuh) \cdot I_{\delta}[\vue -  \vuh]    } \dt. 
\label{A9}
\end{align}
Consequently, keeping in mind that $\Lambda$ and $\delta$ are fixed, we may use the bound \eqref{a6} to conclude 
\begin{align} 
- &\Lambda \int_0^\tau \mathds{1}_{t \in [0,T]} \intO{ I_{\delta}[ \tvue - \vue ] \cdot (\tvue - \wvuh) } \dt + 
\Lambda \int_0^\tau \mathds{1}_{t \in [0,T]}\intO{ I_{\delta}[ \wvuh - \vuh] \cdot  (\tvue - \wvuh) } \dt \br 
&\leq \mathcal{O}(\ep, \tau).
\label{A9a}
\end{align}

\subsubsection{Temperature perturbation}

The nudging term acting on the temperature is represented by the last integral in \eqref{F6}:
\begin{align}
&- \Lambda \int_0^\tau \mathds{1}_{t \in [0,T]} \intO{ \left( \frac{\tvte- \Ov{\vt} }{\ep}-  \wtilde{\mathcal{T}} \right) 
I_\delta\left[ \frac{\tvte - \vte}{\ep} \right] } \dt\br &= - \Lambda \int_0^\tau \mathds{1}_{t \in [0,T]} \intO{ \left( \frac{\tvte - \Ov{\vt} }{\ep}-  \wtilde{\mathcal{T}} \right) 
I_\delta\left[ \frac{\tvte - \Ov{\vt} + \Ov{\vt} - \vte}{\ep} - \wtilde{\mathcal{T}} + \wtilde{\mathcal{T}} \right] } \dt \br 
&= - \Lambda \int_0^\tau \mathds{1}_{t \in [0,T]} \intO{ \left( \frac{\tvte - \Ov{\vt} }{\ep}-  \wtilde{\mathcal{T}} \right) 
I_\delta\left[ \frac{\tvte - \Ov{\vt}}{\ep} - \wtilde{\mathcal{T}} \right] } \dt\br &+ 
\Lambda \int_0^\tau \mathds{1}_{t \in [0,T]} \intO{ \left( \frac{\tvte - \Ov{\vt} }{\ep}-  \wtilde{\mathcal{T}} \right) 
I_\delta\left[ \frac{\vte - \Ov{\vt}}{\ep} - \wtilde{\mathcal{T}} \right] } \dt \br 
&= - \Lambda \int_0^\tau \mathds{1}_{t \in [0,T]} \intO{ \left( \frac{\tvte - \Ov{\vt} }{\ep}-  \wtilde{\mathcal{T}} \right) 
I_\delta\left[ \frac{\tvte - \Ov{\vt}}{\ep} - \wtilde{\mathcal{T}} \right] } \dt \br 
&+ \Lambda \int_0^\tau \mathds{1}_{t \in [0,T]} \intO{ \left( \frac{\tvte - \Ov{\vt} }{\ep}-  \wtilde{\mathcal{T}} \right) 
I_\delta\left[ \frac{\vte - \Ov{\vt}}{\ep} - \mathcal{T} \right] } \dt \br &+ \Lambda \int_0^\tau \mathds{1}_{t \in [0,T]} \intO{ \left( \frac{\tvte - \Ov{\vt} }{\ep}-  \wtilde{\mathcal{T}} \right) 
I_\delta\left[ \mathcal{T} - \wTheta \right] } \dt
\label{A10}
\end{align}

Similarly to Section \ref{vp} above, we have 
\begin{align}
- &\Lambda \int_0^\tau \mathds{1}_{t \in [0,T]} \intO{ \left( \frac{\tvte - \Ov{\vt} }{\ep}-  \wtilde{\mathcal{T}} \right) 
I_\delta\left[ \frac{\tvte - \Ov{\vt}}{\ep} - \wtilde{\mathcal{T}} \right] } \dt \br 
&+ \Lambda \int_0^\tau \mathds{1}_{t \in [0,T]} \intO{ \left( \frac{\tvte - \Ov{\vt} }{\ep}-  \wtilde{\mathcal{T}} \right) 
I_\delta\left[ \frac{\vte - \Ov{\vt}}{\ep} - \mathcal{T} \right] } \dt \leq \mathcal{O}(\ep, \tau). 
\label{A11}
\end{align}
Consequently, it remains to control the integral
\[
\Lambda \int_0^\tau \mathds{1}_{t \in [0,T]} \intO{ \left( \frac{\tvte - \Ov{\vt} }{\ep}-  \wtilde{\mathcal{T}} \right) 
I_\delta\left[ \mathcal{T} - \wTheta \right] } \dt. 
\]

Following \cite{FanFei2025} we denote by $\mathfrak{R}$, $\mathfrak{T}$ the limits 
\begin{align}
\frac{\tvre - \Ov{\vr}}{\ep} &\to \mathfrak{R} \ \mbox{weakly-* in}\ L^\infty(0, T^+; L^{\frac{5}{3}}(\Omega)), \br
\frac{\tvte - \Ov{\vt}}{\ep} &\to \mathfrak{T} \ \mbox{weakly in}\ L^2(0, T^+; W^{1,2}(\Omega)).
\nonumber
\end{align}
Using the estimates \eqref{A9a}, \eqref{A11}, we may proceed exactly as in \cite[Section 5.2]{FanFei2025}
rewriting the relative energy inequality \eqref{F6} in the form 
 \begin{align}
	&\intO{ E_\ep \left(\tvre, \tvte, \tvue \Big| \Ov{\vr} + \ep \wtilde{\mathcal{R}}, \Ov{\vt} + \ep \MTC, \wvuh \right) (\tau, \cdot)}  \br 
	&+ \int_0^\tau \intO{ \Big( \mathbb{S} (\Ov{\vt}, \Ds \tvue) - \mathbb{S} (\Ov{\vt}, \Ds \wvuh) \Big) : \Big( \Ds \tvue - \Ds \wvuh  \Big) } \dt \br &
	+\int_0^\tau \intO{  \left(  \frac{\Ov{\vt} + \ep \MTC}{\tvte^2} \right) \frac{\kappa (\tvte) \Grad \tvte \cdot \Grad \tvte }{\ep^2}  } \dt
	- \int_0^\tau \intO{ \frac{\kappa (\Ov{\vt}) }{\Ov{\vt}} \Grad \mathfrak{T} \cdot \Grad { \MTC} } \dt \br	
	&\leq - \frac{1}{\ep} \int_0^\tau \intO{ \tvre \Big[ (s(\tvre, \tvte) - s( \Ov{\vr} + \ep \wtilde{\mathcal{R}},\Ov{\vt} + \ep \MTC ) \Big] \partial_t \MTC } \dt \br 
	&\quad - \frac{1}{\ep} \int_0^\tau \intO{ \tvre \Big[ s(\tvre, \tvte) - s( \Ov{\vr} + \ep \wtilde{\mathcal{R}}, \Ov{\vt} + \ep \MTC) \Big] \wvuh \cdot \Grad \MTC  } \dt \br
	&\quad  + \frac{1}{\ep} \int_0^\tau \intO{  \tvre \Big[ s(\tvre, \tvte) - s(\Ov{\vr} + \ep \wtilde{\mathcal{R}}, \Ov{\vt} + \ep \MTC) \Big] ( \wvuh - \tvue) \cdot \Grad \MTC } \dt \br
	&\quad + \Lambda \int_0^\tau \mathds{1}_{t \in [0,T]} \intO{ \left( \frac{\tvte - \Ov{\vt} }{\ep}-  \wtilde{\mathcal{T}} \right) 
		I_\delta\left[ \mathcal{T} - \wTheta \right] } \dt \br
	&\quad    +  \int_0^\tau \intO{ \lan \wtilde{\mathcal{R}} - {\mathfrak{R}} \ran \Grad \Big( G + \frac{1}{2}|x_h|^2 \Big)  \cdot \wvuh } \dt  \br 
	&\quad+ C \int_0^\tau \intO{E_\ep \left(\tvre, \tvte, \tvue \Big| \Ov{\vr} + \ep \wtilde{\mathcal{R}}, \Ov{\vt} + \ep \MTC, \wvuh \right) } \dt   + \mathcal{O}(\ep, \tau),
	\label{A16}
\end{align}
cf. \cite[formula (5.25)]{FanFei2025}.

Recall that $\wTheta$ satisfies \eqref{A5}, specifically
\begin{align} 
\Ov{\vr} c_p(\Ov{\vr}, \Ov{\vt} ) \Big[ \partial_t \wTheta + \wvuh \cdot \Gradh \wTheta \Big] &- 
\Ov{\vr} \Ov{\vt} \alpha (\Ov{\vr}, \Ov{\vt} ) \wvuh \cdot \Gradh \Big( G + \frac{1}{2}|\vc x_h|^2 \Big) \br &= 
\kappa (\Ov{\vt}) \Del \wTheta + \Ov{\vt}  \alpha (\Ov{\vr}, \Ov{\vt} )  \frac{\partial p  (\Ov{\vr}, \Ov{\vt} ) } 
{\partial \vt} \partial_t \avintO{ \wTheta } - \Ov{\vt} \Lambda I_\delta \left[ \wTheta - \mathcal{T} \right]\mathds{1}_{t \in [0,T]}. 
\label{A17}
\end{align}
Next, observe that 
\begin{align}
	\frac{1}{\ep} &\int_0^\tau \intO{  \tvre \Big[ s(\tvre,\tvte) - s(\Ov{\vr} + \ep \wtilde{\mathcal{R}}, \Ov{\vt} + \ep \MTC) \Big] ( \wvuh - \tvue) \cdot \Grad {\MTC} } \dt \br
	&\aleq \int_0^\tau \intO{E_\ep \left(\tvre, \tvte, \tvue \Big| \Ov{\vr} + \ep \wtilde{\mathcal{R}}, \Ov{\vt} + \ep \MTC, \wvuh \right) } \dt,
	\nonumber
\end{align}
and 
\begin{align}
\Lambda &\int_0^\tau \mathds{1}_{t \in [0,T]} \intO{ \left( \frac{\tvte - \Ov{\vt} }{\ep}-  \wtilde{\mathcal{T}} \right) 
	I_\delta\left[ \mathcal{T} - \wTheta \right] } \dt \br &\leq \Lambda \int_0^\tau \mathds{1}_{t \in [0,T]} \intO{ \left( \mathfrak{T}-  \wtilde{\mathcal{T}} \right) 
	I_\delta\left[ \mathcal{T} - \wTheta \right] } \dt + \mathcal{O}(\ep, \tau).
\nonumber	
\end{align}
Consequently, inequality \eqref{A16} reduces to 
 \begin{align}
	&\intO{ E_\ep \left(\tvre, \tvte, \tvue \Big| \Ov{\vr} + \ep \wtilde{\mathcal{R}}, \Ov{\vt} + \ep \MTC, \wvuh \right) (\tau, \cdot)}  \br 
	&+ \int_0^\tau \intO{ \Big( \mathbb{S} (\Ov{\vt}, \Ds \tvue) - \mathbb{S} (\Ov{\vt}, \Ds \wvuh) \Big) : \Big( \Ds \tvue - \Ds \wvuh  \Big) } \dt \br &
	+\int_0^\tau \intO{  \left(  \frac{\Ov{\vt} + \ep \MTC}{\tvte^2} \right) \frac{\kappa (\tvte) \Grad \tvte \cdot \Grad \tvte }{\ep^2}  } \dt
	- \int_0^\tau \intO{ \frac{\kappa (\Ov{\vt}) }{\Ov{\vt}} \Grad \mathfrak{T} \cdot \Grad { \MTC} } \dt \br	
	&\leq - \int_0^\tau \intO{  \Ov{\vr} \left(\frac{\partial s(\Ov{\vr}, \Ov{\vt})}{\partial \vr}(\mathfrak{R} - \wtilde{\mathcal{R}}) + \frac{\partial s(\Ov{\vr}, \Ov{\vt})}{\partial \vt}(\mathfrak{T} - \MTC) \right) {\frac{\Ov{\vt} \alpha (\Ov{\vr}, \Ov{\vt} )}{ c_p (\Ov{\vr}, \Ov{\vt} )} } \Grad \Big( G + | x_h|^2 \Big) \cdot \wvuh                       } \dt \br
	&\quad - { \int_0^\tau \intO{  \Ov{\vr} \left(\frac{\partial s(\Ov{\vr}, \Ov{\vt})}{\partial \vr}(\mathfrak{R} - \wtilde{\mathcal{R}}) + \frac{\partial s(\Ov{\vr}, \Ov{\vt})}{\partial \vt}(\mathfrak{T} - \MTC) \right)  \frac{\kappa(\Ov{\vt})}{\Ov{\vr} c_p (\Ov{\vr}, \Ov{\vt} )} \Del \MTC                  } \dt} \br	
	&\quad { -  \int_0^\tau \intO{  \Ov{\vr} \left(\frac{\partial s(\Ov{\vr}, \Ov{\vt})}{\partial \vr}(\mathfrak{R} - \wtilde{\mathcal{R}}) + \frac{\partial s(\Ov{\vr}, \Ov{\vt})}{\partial \vt}(\mathfrak{T} - \MTC) \right)  \frac{1}{\Ov{\vr} c_p (\Ov{\vr}, \Ov{\vt} )} \xi(t)                  } \dt} \br
	&\quad { +  \int_0^\tau \mathds{1}_{t \in [0,T] }\intO{  \Ov{\vr} \left(\frac{\partial s(\Ov{\vr}, \Ov{\vt})}{\partial \vr}(\mathfrak{R} - \wtilde{\mathcal{R}}) + \frac{\partial s(\Ov{\vr}, \Ov{\vt})}{\partial \vt}(\mathfrak{T} - \MTC) \right)  \frac{\Ov{\vt} \Lambda }{\Ov{\vr} c_p (\Ov{\vr}, \Ov{\vt} )} I_\delta[ \wTheta - \mathcal{T}]  } \dt} \br
	&\quad + \Lambda \int_0^\tau \mathds{1}_{t \in [0,T]} \intO{ \left( \mathfrak{T} -  \wtilde{\mathcal{T}} \right) 
		I_\delta\left[ \mathcal{T} - \wTheta \right] } \dt \br
	&\quad    +  \int_0^\tau \intO{ \lan \wtilde{\mathcal{R}} - {\mathfrak{R}} \ran \Grad \Big( G + \frac{1}{2}| x_h|^2 \Big)  \cdot \wvuh } \dt  \br 
	&\quad+ C \int_0^\tau \intO{E_\ep \left(\tvre, \tvte, \tvue \Big| \Ov{\vr} + \ep \wtilde{\mathcal{R}}, \Ov{\vt} + \ep \MTC, \wvuh \right) } \dt   + \mathcal{O}(\ep, \tau),
	\label{A18}
\end{align}
where we have denoted 
\[
\xi = \Ov{\vt}  \alpha (\Ov{\vr}, \Ov{\vt} )  \frac{\partial p  (\Ov{\vr}, \Ov{\vt} ) } 
{\partial \vt} \partial_t \avintO{ \wTheta }.
\]

Since both $\mathfrak{R}$, $\mathfrak{T}$ and $\wtilde{\mathcal{R}}$, $\wTheta$ satisfy the Boussinesq relation \eqref{aa2}, 
we get 
\begin{equation} \label{A19}
		\frac{\partial p(\Ov{\vr}, \Ov{\vt})}{\partial \vr} \Grad (\wtilde{\mathcal{R}} - \mathfrak{R}) + \frac{\partial p(\Ov{\vr}, \Ov{\vt})}{\partial \vt} \Grad (\wTheta - \mathfrak{T}) = 0
	,\ \intO{ ( \wtilde{\mathcal{R}} - \mathfrak{R} ) } = 0.
\end{equation}	

Next, as 
\[
\intO{ ( \wtilde{\mathcal{R}} - \mathfrak{R} ) } = 0, 
\]
we get 
\begin{align}
	-  \int_0^\tau \intO{  \Ov{\vr} \left(\frac{\partial s(\Ov{\vr}, \Ov{\vt})}{\partial \vr}(\mathfrak{R} - \wtilde{\mathcal{R}}) + \frac{\partial s(\Ov{\vr}, \Ov{\vt})}{\partial \vt}(\mathfrak{T} - \MTC) \right)  \frac{1}{\Ov{\vr} c_p (\Ov{\vr}, \Ov{\vt} )} \xi(t)                  } \dt \nonumber \br 
	=- \int_0^\tau \intO{  \frac{\partial s(\Ov{\vr}, \Ov{\vt})}{\partial \vt} \left[
		\frac{\partial p(\Ov{\vr}, \Ov{\vt} )}{\partial \vr} \left( \frac{\partial p(\Ov{\vr}, \Ov{\vt} )}{\partial \vt} \right)^{-1}(\mathfrak{R} - \wtilde{\mathcal{R}}) + (\mathfrak{T} - \MTC) \right]  \frac{1}{ c_p (\Ov{\vr}, \Ov{\vt} )} \xi(t)                  } \dt,
	\label{A20}
\end{align}
where, thanks to \eqref{A19}, the quantity 
\[
\left[
\frac{\partial p(\Ov{\vr}, \Ov{\vt} )}{\partial \vr} \left( \frac{\partial p(\Ov{\vr}, \Ov{\vt} )}{\partial \vt} \right)^{-1}(\mathfrak{R} - \wtilde{\mathcal{R}}) + (\mathfrak{T} - \MTC) \right]
\ \mbox{is independent of}\ x \in \Omega. 
\]

Similarly, 
\begin{align} 
 &- { \int_0^\tau \intO{  \Ov{\vr} \left(\frac{\partial s(\Ov{\vr}, \Ov{\vt})}{\partial \vr}(\mathfrak{R} - \wtilde{\mathcal{R}}) + \frac{\partial s(\Ov{\vr}, \Ov{\vt})}{\partial \vt}(\mathfrak{T} - \MTC) \right)  \frac{\kappa(\Ov{\vt})}{\Ov{\vr} c_p (\Ov{\vr}, \Ov{\vt} )} \Del \MTC } \dt} \br	
&\quad +  \int_0^\tau \mathds{1}_{t \in [0,T] }\intO{  \Ov{\vr} \left(\frac{\partial s(\Ov{\vr}, \Ov{\vt})}{\partial \vr}(\mathfrak{R} - \wtilde{\mathcal{R}}) + \frac{\partial s(\Ov{\vr}, \Ov{\vt})}{\partial \vt}(\mathfrak{T} - \MTC) \right)  \frac{\Ov{\vt} \Lambda }{\Ov{\vr} c_p (\Ov{\vr}, \Ov{\vt} )} I_\delta[ \wTheta - \mathcal{T}]  } \dt \br 
&\quad = 	-	\int_0^\tau \intO{  \frac{\partial s(\Ov{\vr}, \Ov{\vt})}{\partial \vr}\left[ (\mathfrak{R} - \wtilde{\mathcal{R}}) +  \frac{\partial p(\Ov{\vr}, \Ov{\vt}) }{\partial \vt}
		\left(  \frac{\partial p(\Ov{\vr}, \Ov{\vt}) }{\partial \vr} \right)^{-1} (\mathfrak{T} - \MTC)  \right]  \frac{\kappa(\Ov{\vt})}{ c_p (\Ov{\vr}, \Ov{\vt} )} \Del \MTC                  } \dt \br
&\quad 		
+	\int_0^\tau \mathds{1}_{t \in [0,T] }\intO{  \frac{\partial s(\Ov{\vr}, \Ov{\vt})}{\partial \vr}\left[ (\mathfrak{R} - \wtilde{\mathcal{R}}) +  \frac{\partial p(\Ov{\vr}, \Ov{\vt}) }{\partial \vt}
		\left(  \frac{\partial p(\Ov{\vr}, \Ov{\vt}) }{\partial \vr} \right)^{-1} (\mathfrak{T} - \MTC)  \right]  \frac{ \Ov{\vt} \Lambda
		}{ c_p (\Ov{\vr}, \Ov{\vt} )} I_\delta[ \wTheta - \mathcal{T}]                  } \dt		\br
&\quad	+ \int_0^\tau \intO{ \left(  \frac{\partial s(\Ov{\vr}, \Ov{\vt})}{\partial \vr}\frac{\partial p(\Ov{\vr}, \Ov{\vt}) }{\partial \vt}
		\left(  \frac{\partial p(\Ov{\vr}, \Ov{\vt}) }{\partial \vr} \right)^{-1} (\mathfrak{T} - \MTC) - \frac{\partial s(\Ov{\vr}, \Ov{\vt})}{\partial \vt} (\mathfrak{T} - \MTC)\right) \frac{\kappa(\Ov{\vt})} { c_p (\Ov{\vr}, \Ov{\vt} )} \Del \MTC  } \dt \br 
&\quad - 	\int_0^\tau \mathds{1}_{t \in [0,T] } \times \br 
\quad &\times \intO{ \left(  \frac{\partial s(\Ov{\vr}, \Ov{\vt})}{\partial \vr}\frac{\partial p(\Ov{\vr}, \Ov{\vt}) }{\partial \vt}
		\left(  \frac{\partial p(\Ov{\vr}, \Ov{\vt}) }{\partial \vr} \right)^{-1} (\mathfrak{T} - \MTC) - \frac{\partial s(\Ov{\vr}, \Ov{\vt})}{\partial \vt} (\mathfrak{T} - \MTC)\right) \frac{ \Ov{\vt} \Lambda
		}{ c_p (\Ov{\vr}, \Ov{\vt} )} I_\delta[ \wTheta - \mathcal{T}]                } \dt,
 \label{A21}	
\end{align}
where, similarly to the above, the quantity 
\begin{equation} \label{A22}
\left[ (\mathfrak{R} - \wtilde{\mathcal{R}}) +  \frac{\partial p(\Ov{\vr}, \Ov{\vt}) }{\partial \vt}
		\left(  \frac{\partial p(\Ov{\vr}, \Ov{\vt}) }{\partial \vr} \right)^{-1} (\mathfrak{T} - \MTC)  \right]
		\ \mbox{is independent of}\ x \in \Omega.
\end{equation}

Thanks to \eqref{A22}, the value of the first two integrals on the right--hand side of \eqref{A21} reduces to the 
integral mean of the quantity
\[
\frac{\kappa(\Ov{\vt})}{ c_p (\Ov{\vr}, \Ov{\vt} )} \Del \MTC    - \frac{ \Ov{\vt} \Lambda
}{ c_p (\Ov{\vr}, \Ov{\vt} )} I_\delta[ \wTheta - \mathcal{T}]                     
\]
that can be computed integrating equation \eqref{A17},
\begin{align}
&\avintO{ \left( \frac{\kappa(\Ov{\vt})}{ c_p (\Ov{\vr}, \Ov{\vt} )} \Del \MTC    - \frac{ \Ov{\vt} \Lambda
}{ c_p (\Ov{\vr}, \Ov{\vt} )} I_\delta[ \wTheta - \mathcal{T}]\right) }	= 
\left( \Ov{\vr} - \Ov{\vt} \frac{\alpha (\Ov{\vr}, \Ov{\vt})}{c_p  (\Ov{\vr}, \Ov{\vt}) }   
\frac{\partial p  (\Ov{\vr}, \Ov{\vt}) }{\partial \vt} \right) \avintO{ \partial_t \wTheta } \br 
&\quad = \left( \Ov{\vr} - \Ov{\vt} \frac{\alpha (\Ov{\vr}, \Ov{\vt})}{c_p  (\Ov{\vr}, \Ov{\vt}) }   
\frac{\partial p  (\Ov{\vr}, \Ov{\vt}) }{\partial \vt} \right) \left( \Ov{\vt} {\alpha (\Ov{\vr}, \Ov{\vt})} 
\frac{\partial p  (\Ov{\vr}, \Ov{\vt}) }{\partial \vt} \right)^{-1} \xi.
\nonumber
\end{align}	
The resulting expression can be added to \eqref{A20}. 

The last two integrals in \eqref{A21} can be computed exactly as in \cite[Section 5.2]{FanFei2025}:
\begin{align} 
	\int_0^\tau & \intO{\left[ \left( \frac{\partial s(\Ov{\vr}, \Ov{\vt})}{\partial \vr}\frac{\partial p(\Ov{\vr}, \Ov{\vt}) }{\partial \vt}
		\left(  \frac{\partial p(\Ov{\vr}, \Ov{\vt}) }{\partial \vr} \right)^{-1} - \frac{\partial s(\Ov{\vr}, \Ov{\vt})}{\partial \vt} \right)(\mathfrak{T} - \MTC) \right]  \frac{\kappa(\Ov{\vt})} { c_p (\Ov{\vr}, \Ov{\vt} )} \Del \MTC  } \dt \br 
&= \frac{\kappa(\Ov{\vt})} { \Ov{\vt} }\int_0^\tau \intO{ (\mathfrak{T} - \wTheta) \Del \wTheta       } \dt, 		
\nonumber
\end{align}
and 
\begin{align}
 - 	&\int_0^\tau \mathds{1}_{t \in [0,T] }  \intO{ \left[ \left( \frac{\partial s(\Ov{\vr}, \Ov{\vt})}{\partial \vr}\frac{\partial p(\Ov{\vr}, \Ov{\vt}) }{\partial \vt}
	\left(  \frac{\partial p(\Ov{\vr}, \Ov{\vt}) }{\partial \vr} \right)^{-1} - \frac{\partial s(\Ov{\vr}, \Ov{\vt})}{\partial \vt} \right)(\mathfrak{T} - \MTC) \right] \frac{ \Ov{\vt} \Lambda
	}{ c_p (\Ov{\vr}, \Ov{\vt} )} I_\delta[ \wTheta - \mathcal{T}]                } \dt \br 
&= - \Lambda \int_0^\tau \mathds{1}_{t \in [0,T] } \intO{ (\mathfrak{T} - \wTheta) I_\delta [ \mathcal{T} - \wTheta]    }. 
\label{A23}
\end{align}

Now, the integral \eqref{A23} cancels out with its counterpart in \eqref{A18}. 
Thus repeating step by step the arguments of \cite[Section 5.2]{FanFei2025} we obtain the desired conclusion 
\begin{align}
	&\intO{ E_\ep \left(\tvre, \tvte, \tvue \Big| \Ov{\vr} + \ep \wtilde{\mathcal{R}}, \Ov{\vt} + \ep \MTC, \wvuh \right) (\tau, \cdot) }  \br 
	&+ \int_0^\tau \intO{ \Big( \mathbb{S} (\Ov{\vt}, \Ds \tvue) - \mathbb{S} (\Ov{\vt}, \Ds \wvuh) \Big) : \Big( \Ds \tvue - \Ds \wvuh  \Big) } \dt \br  &+ \int_0^\tau \intO{ \frac{\kappa (\Ov{\vt} ) }{\Ov{\vt}} \left|\Grad 
		\left(	\frac{\tvte - \Ov{\vt}}{\ep} \right) - \Grad \MTC \right|^2 } \dt \br
	&\aleq \int_0^\tau \intO{E_\ep \left(\tvre, \tvte, \tvue \Big|\Ov{\vr} + \ep \wtilde{\mathcal{R}}, \Ov{\vt} + \ep \MTC, \wvuh \right) } \dt   + \mathcal{O}(\ep, \tau).
	\label{A24}
\end{align}
Thus the standard application of the Gr\" onwall lemma yields \eqref{F2}.

\section*{Compliance with Ethical Standards}\label{conflicts}

\smallskip
\par\noindent 

{\bf Conflict of Interest}. The authors declare that they have no conflict of interest.

\smallskip
\par\noindent
{\bf Data Availability}. Data sharing is not applicable to this article as no datasets were generated or analysed during the current study.

\def\cprime{$'$} \def\ocirc#1{\ifmmode\setbox0=\hbox{$#1$}\dimen0=\ht0
	\advance\dimen0 by1pt\rlap{\hbox to\wd0{\hss\raise\dimen0
			\hbox{\hskip.2em$\scriptscriptstyle\circ$}\hss}}#1\else {\accent"17 #1}\fi}


\end{document}